\documentclass[12pt]{article}
\usepackage{times}
\usepackage{amsmath}
\usepackage{amssymb}
\usepackage{amsthm}
\usepackage{amsmath}
\usepackage{amsfonts}
\usepackage{amssymb}
\usepackage[all]{xy}
\usepackage{xcolor}
\usepackage[utf8]{inputenc}

\newtheorem{theorem}{Theorem}[section]
\newtheorem{lemma}[theorem]{Lemma}
\newtheorem{corollary}[theorem]{Corollary}
\newtheorem{proposition}[theorem]{Proposition}

\theoremstyle{definition}
\newtheorem{observation}[theorem]{Remark}
\newtheorem{definition}[theorem]{Definition}
\def \fim {{\hfill $\blacksquare$}}

\def\Z{\mathbb Z}

\def\0{\underline 0}

\begin{document}

\title {\bf Brasselet number and function-germs with a one-dimensional critical set   \footnote{Research partially supported by  FAPESP - Brazil, Grant 2015/25191-9 and 2017/18543-1. \newline $\quad$ {\it
Key-words: Brasselet number, Euler obstruction, stratified Morse critical points}
\newline   {\it } }}

\vspace{1cm}
\author{Hellen Santana }

\date{\bf }
\maketitle

\begin{abstract}
The Brasselet number of a function $f$ with nonisolated singularities describes numerically the topological information of its generalized Milnor fibre. In this work, using the Brasselet number, we present several formulas for germs $f:(X,0)\rightarrow(\mathbb{C},0)$ and $g:(X,0)\rightarrow(\mathbb{C},0)$ in the case where $g$ has a one-dimensional critical locus. We also give applications when $f$ has isolated singularities and when it is a generic linear form.

\end{abstract}

\section*{Introduction}

\hspace{0,5cm} Let $f:(\mathbb{C}^n,0)\rightarrow(\mathbb{C},0)$ be an analytic function defined in a neighborhood of the origin and $\Sigma f$ the critical locus of $f.$ Milnor studied the set  $f^{-1}(\delta)\cap B_{\epsilon}$, denoted by $F_{f,0}$ and later called Milnor fiber, where $\delta$ is a regular value of $f, 0<|\delta|\ll\epsilon\ll1.$ In \cite{Milnor}, Milnor proved that, if $f$ has an isolated singularity, $F_{f,0}$ has the homotopy type of a wedge of $\mu(f)$ spheres of dimension $n-1,$ where $\mu(f)$ is the Milnor number of $f.$ This number also gives an important geometric information associated to the function $f,$ which is the number of Morse points in a Morsification of $f$ in a neighborhood of the origin. 

In \cite{Hamm}, Hamm generalized Milnor's results for complete intersections with isolated singularity $F=(f_1,\ldots,f_k):(\mathbb{C}^n,0)\rightarrow(\mathbb{C}^k,0), 1<k<n,$ proving that the Milnor fiber $F^{-1}(\delta)\cap B_{\epsilon}, 0<|\delta|\ll\epsilon\ll1,$ has the homotopy type of a wedge of $\mu(F)$ spheres of dimension $n-k.$ In this context, Lê \cite{L} and Greuel \cite{Greuel} proved that $\mu(F)+\mu(F^{\prime})=\dim_{\mathbb{C}}(\frac{\mathcal{O}_{\mathbb{C}^n,0}}{I}),$ where $F^{\prime}:(\mathbb{C}^n,0)\rightarrow(\mathbb{C}^{k-1},0)$ is the map with components $f_1,\ldots,f_{k-1}$ and $I$ is the ideal generated by $f_1,\ldots,f_{k-1}$ and the $(k\times k)-$ minors $\frac{\partial(f_1,\ldots,f_{k})}{\partial(x_{i_1},\ldots,x_{i_k})}.$ Notice that the number $\dim_{\mathbb{C}}(\frac{\mathcal{O}_{\mathbb{C}^n,0}}{I})$ is the number of critical points of a Morsification of $f_k$ appearing on the Milnor fibre of $F^{\prime}.$ 

If $f$ is defined over a complex analytic space $X$ and $f$ has an isolated singularity at the origin, a generalization for the Milnor number is the Euler obstruction of the function $f$, introduced in \cite{BMPS}, by Brasselet, Massey, Parameswaran and Seade. In \cite{STV}, Seade, Tib\u{a}r and Verjovsky proved that, up to sign, this number is the number of Morse critical points of a stratified Morsification of $f$ appearing in the regular part of $X$ in a neighborhood of the origin.  

In a more general context, if $f$ is defined over a complex analytic germ $(X,0)$ equipped with a good stratification $\mathcal{V}$ relative to $f$ (see Definition \ref{good stratification}) and the function $f$ does not have isolated singularity at the origin, a way to describe the generalized Milnor fiber $X\cap f^{-1}(\delta)\cap B_{\epsilon}$ is using the Brasselet number of $f$ at the origin, $B_{f,X}(0),$ introduced by Dutertre and Grulha in \cite{DG}. In that paper, the authors presented a Lê-Greuel type formula for the Brasselet number: if $g:X\rightarrow\mathbb{C}$ is prepolar with respect to $\mathcal{V}$ at the origin (see Definition \ref{definition prepolar}) and $0<|\delta|\ll\epsilon\ll1,$ then
\begin{center}
$B_{f,X}(0)-B_{f,X^g}(0)=(-1)^{d-1}n_q,$
\end{center} 
where $n_q$ is the number of Morse critical points of a partial Morsification of $g|_{X\cap f^{-1}(\delta)\cap B_{\epsilon}}$ appearing in the regular part of $X$, and $X^g=X\cap \{g=0\}.$

They also proved several results about the topology of functions with isolated singularity defined over an analytic complex Whitney stratified variety $X.$ If $X$ is equidimensional, let $f,g:X\rightarrow\mathbb{C}$ be analytic functions with isolated singularity at the origin, such that $g$ is prepolar with respect to the good stratification induced by $f$ at the origin (see Definition \ref{induced stratification}) and $f$ is prepolar with respect to the good stratification induced by $g$ at the origin, then $B_{f,X^g}(0)=B_{g,X^f}(0),$ where $X^f=X\cap\{f=0\}.$ Also, if $n_q$ is the number of Morse critical points of a Morsification of $g|_{X\cap f^{-1}(\delta)\cap B_{\epsilon}}$ appearing in the regular part of $X$ and $m_q$ is the number of Morse critical points of a Morsification of $f|_{X\cap g^{-1}(\delta)\cap B_{\epsilon}}$ appearing in the regular part of $X$, for $0<|\delta|\ll\epsilon\ll 1,$ then\begin{center}
$B_{f,X}(0)-B_{g,X}(0)=(-1)^{d-1}(n_q-m_q).$
\end{center}

An interesting consequence of this last statement is a way to compare the local Euler obstruction $Eu_{X^g}(0)$ and the Brasselet number $B_{g,X\cap H}(0),$ given by the equality $Eu_{X^g}(0)=B_{g,X\cap H}(0),$ where $H$ is a generic hyperplane passing through the origin.

In this work, we start considering, in Section 3, two function-germs $f,g:X\rightarrow\mathbb{C}$ and a good stratification $\mathcal{V}$ of $X$ relative to $f.$ We suppose that the critical locus $\Sigma_{\mathcal{V}}g$ of $g$ is one-dimensional, that $\Sigma_{\mathcal{V}}g\cap\{f=0\}=\{0\}$ and we denote by $\mathcal{V}^f$ the collection of strata of $\mathcal{V}$ contained in $\{f=0\}$ and by $V_1,\ldots, V_q$ the strata of $\mathcal{V}$ not contained in $\{f=0\}.$  Then, we prove (Lemma \ref{first stratification lemma}) that the refinement \begin{eqnarray}
\mathcal{V}^{\prime}=\Big\{V_{i}\setminus \Sigma_{\mathcal{V}} g, V_{i}\cap \Sigma_{\mathcal{V}} g, i\in\{1,\ldots,q\}\Big\}\cup\mathcal{V}^f\label{good stratification}
\end{eqnarray} 
is a good stratification of $X$ relative to $f$ and  $\mathcal{V}^{\prime\{g=0\}}$ is a good stratification of \linebreak $X\cap\{g=0\}$ relative to $f|_{X\cap\{g=0\}},$ where \begin{center}$\mathcal{V}^{\prime\{g=0\}}=\Big\{V_{i}\cap\{g=0\} \setminus \Sigma_{\mathcal{V}} g, V_{i}\cap \Sigma_{\mathcal{V}} g, i\in\{1,\ldots,q\}\Big\}\cup\Big(\mathcal{V}^f\cap\{g=0\}\Big),$\end{center}
\noindent and $\mathcal{V}^f\cap\{g=0\}$ denotes the collection of strata of type $V^f\cap\{g=0\},$ with $V^f\in\mathcal{V}^f.$
We write $\Sigma_{\mathcal{V}} g$ as a union of irreducible components (branches) \linebreak$\Sigma_{\mathcal{V}} g=b_1\cup\ldots\cup b_r,$ where $b_j\subseteq V_{i_j},$ for some $i_j\in\{1,\ldots,q\}$ and we take a regular value $\delta$ of $f, 0<|\delta|\ll1,$ and, for each $j\in\{1,\ldots,r\},$ we set $f^{-1}(\delta)\cap b_j=\{x_{i_1},\ldots,x_{i_{k(j)}}\}.$ So, in this case, the local degree $m_{f,b_j}$ of $f|_{b_j}$ is $k.$ Let $\epsilon$ be sufficiently small such that the local Euler obstruction of $X$ and $X^g$ are constant on $b_j\cap B_{\epsilon}$. In this case, we denote by $Eu_X(b_j)$ (respectively, $Eu_{X^g}(b_j)$) the local Euler obstruction of $X$ (respectively, $X^g$) at a point of $b_j\cap B_{\epsilon}.$ If $g$ is tractable at the origin with respect to $\mathcal{V}$ relative to $f$ (see Definition \ref{definition tractable}) and $0<|\delta|\ll\epsilon\ll1,$ we prove (Theorem \ref{Generalization 6.4 of DG}) that
\begin{center}
$B_{f,X}(0)-B_{f,X^g}(0)-\sum_{j=1}^{r}m_{f,b_j}(Eu_X(b_j)-Eu_{X^g}(b_j))=(-1)^{d-1}m,$
\end{center}
where $m$ is the number of stratified Morse critical points of a Morsification of \linebreak$g:X\cap f^{-1}(\delta)\cap B_{\epsilon}\rightarrow\mathbb{C}$ appearing on $X_{reg}\cap f^{-1}(\delta)\cap \{g\neq 0\}\cap B_{\epsilon}.$

We conclude that, in the case where $g$ is not prepolar with respect to $\mathcal{V}$ relative to $f$, the Lê-Greuel type formula for the Brasselet number presents a type of defect. More precisely, this formula shows us that the number of Morse critical points $m$ on the regular part of $X$ does not contain all the topological information given by the difference $B_{f,X}(0)-B_{f,X^g}(0).$

In Section 4, we suppose that $f$ has an isolated singularity at the origin and we consider a Whitney stratification $\mathcal{W}$ of $X.$ Let $\mathcal{V}$ be the good stratification of $X$ induced by $f$ and let us suppose that $g$ is tractable at the origin with respect to $\mathcal{V}$ relative to $f.$ We prove (Lemma \ref{second stratification lemma}) that the refinement $\mathcal{V}^{\prime\prime}$ of $\mathcal{V},$\begin{eqnarray}
\mathcal{V}^{\prime\prime}=\Big\{V_i\setminus\{g=0\}, V_i\cap\{g=0\}\setminus\Sigma_{\mathcal{W}}g, V_i\cap\Sigma_{\mathcal{W}}g, V_i\in\mathcal{V}\Big\}\cup\{0\}
\end{eqnarray} is a good stratification of $X$ relative to $g$ such that 
$\mathcal{V}^{\prime\prime\{f=0\}}=\{V_i^{\prime\prime}\cap\{f=0\},V_i^{\prime\prime}\in\mathcal{V}^{\prime\prime} \}$ is a good stratification of $X^f$ relative to $g|_{X^f}.$ 

Using this stratification, we prove (Corollary \ref{Brasselet number 1}) that, \begin{center}
$B_{g,X^f}(0)=B_{f,X^g}(0)-\sum_{j=1}^{r}m_{f,b_j}(Eu_{X^g}(b_j)-B_{g,X\cap \{f=\delta\}}(b_j)),$
\end{center}

\noindent where $B_{g,X\cap \{f=\delta\}}(b_j)$ is the Brasselet number of $g$ at a point $x_{j_s}\in b_j\cap\{f=\delta\}.$

As a consequence of this result, we obtain a way to compare the local Euler obstruction $Eu_{X^g}(0)$ and the Brasselet number $B_{g,X\cap H}(0)$ in the case where $g$ has one-dimensional critical locus. Let $l$ be a generic linear form over $\mathbb{C}^n$ and $H=l^{-1}(0).$ We prove (Corollary \ref{generalization 6.6}) that: \begin{center}
$B_{g,X\cap H}(0)=Eu_{X^g}(0)-\sum_{j=1}^{r}m_{b_j}(Eu_{X^g}(b_j)-B_{g,X\cap l^{-1}(\delta)}(b_j)),$ \end{center} 
\noindent where $m_{b_j}$ is the multiplicity of the branch $b_j$ at the origin.
In this same setting, we also prove (Corollary \ref{generalization 6.5}) that \begin{center}
$B_{g,X}(0)-B_{f,X}(0)=(-1)^{d-1}(n_q-m_q)-\sum_{j=1}^{r}m_{f,b_j}(Eu_X(b_j)-B_{g,X\cap\{f=\delta\}}(b_j)).$
\end{center}

In Section 1, we present definitions and results about objects we will need to develop this work, like the local Euler obstruction and the Euler obstruction of a function. Section 2 is dedicated to the central element of this paper, the Brasselet number.

\section{Local Euler obstruction and Euler obstruction of a function}

\hspace{0,5cm} In this section, we will see the definition of the local Euler obstruction, a singular invariant defined by MacPherson and used as one of the main tools in his proof of the Deligne-Grothendieck conjecture about the existence and uniqueness of Chern classes for singular varities. 

Let $(X,0)\subset(\mathbb{C}^n,0)$ be an equidimensional reduced complex analytic germ of dimension $d$ in a open set $U\subset\mathbb{C}^n.$ Consider a complex analytic Whitney stratification $\mathcal{V}=\{V_{\lambda}\}$ of $U$ adapted to $X$ such that $\{0\}$ is a stratum. We choose a small representative of $(X,0),$ denoted by $X,$ such that $0$ belongs to the closure of all strata. We write $X=\cup_{i=0}^{q} V_i,$ where $V_0=\{0\}$ and $V_q=X_{reg},$ where $X_{reg}$ is the regular part of $X.$ We suppose that $V_0,V_1,\ldots,V_{q-1}$ are connected and that the analytic sets $\overline{V_0},\overline{V_1},\ldots,\overline{V_q}$ are reduced. We write $d_i=dim(V_i), \ i\in\{1,\ldots,q\}.$ Note that $d_q=d.$ 

Let $G(d,N)$ be the Grassmannian manifold, $x\in X_{reg}$ and consider the Gauss map $\phi: X_{reg}\rightarrow U\times G(d,N)$ given by $x\mapsto(x,T_x(X_{reg})).$ 

\begin{definition}
The closure of the image of the Gauss map $\phi$ in $U\times G(d,N)$, denoted by $\tilde{X}$, is called \textbf{Nash modification} of $X$. It is a complex analytic space endowed with an analytic projection map $\nu:\tilde{X}\rightarrow X.$
\end{definition}

Consider the extension of the tautological bundle $\mathcal{T}$ over $U\times G(d,N).$ Since \linebreak$\tilde{X}\subset U\times G(d,N)$, we consider $\tilde{T}$ the restriction of $\mathcal{T}$ to $\tilde{X},$ called the \textbf{Nash bundle}, and $\pi:\tilde{T}\rightarrow\tilde{X}$ the projection of this bundle.

In this context, denoting by $\varphi$ the natural projection of $U\times G(d,N)$ at $U,$ we have the following diagram:

$$\xymatrix{
\tilde{T} \ar[d]_{\pi}\ar[r] & \mathcal{T}\ar[d] \\ 
\tilde{X}\ar[d]_{\nu}\ar[r] & U\times G(d,N)\ar[d]^{\varphi} \\ 
X\ar[r] & U\subseteq\mathbb{C}^N \\}  $$

Considering $\vert\vert z\vert\vert=\sqrt{z_1\overline{z_1}+\cdots+z_N\overline{z_N}}$, the $1$-differential form $w=d\vert\vert z\vert\vert^2$ over $\mathbb{C}^N$ defines a section in $T^{*}\mathbb{C}^N$ and its pullback $\varphi^{*}w$ is a $1$- form over $U\times G(d,N).$ Denote by $\tilde{w}$ the restriction of $\varphi^{*}w$ over $\tilde{X}$, which is a section of the dual bundle $\tilde{T}^{*}.$

Choose $\epsilon$ small enough for $\tilde{w}$ be a non zero section over $\nu^{-1}(z), 0<\vert\vert z \vert\vert\leqslant\epsilon,$ let $B_{\epsilon}$ be the closed ball with center at the origin with radius $\epsilon$ and denote by:

\begin{enumerate}

\item $Obs(\tilde{T}^{*},\tilde{w})\in\mathbb{H}^{2d}(\nu^{-1}(B_{\epsilon}),\nu^{-1}(S_{\epsilon}),\Z)$  the obstruction for extending $\tilde{w}$ from $\nu^{-1}(S_{\epsilon})$ to $\nu^{-1}(B_{\epsilon});$

\item $O_{\nu^{-1}(B_{\epsilon}),\nu^{-1}(S_{\epsilon})}$ the fundamental class in $\mathbb{H}_{2d}(\nu^{-1}(B_{\epsilon}),\nu^{-1}(S_{\epsilon}),\Z).$ 
\end{enumerate}

\begin{definition}
The \textbf{local Euler obstruction} of $X$ at $0, \ Eu_X(0),$ is given by the evaluation $$Eu_X(0)=\langle Obs(\tilde{T}^{*},\tilde{w}),O_{\nu^{-1}(B_{\epsilon}),\nu^{-1}(S_{\epsilon})}\rangle.$$
\end{definition}

In \cite{BLS}, Brasselet, Lê and Seade proved a formula to make the calculation of the Euler obstruction easier.

\begin{theorem}(Theorem 3.1 of \cite{BLS})
Let $(X,0)$ and $\mathcal{V}$ be given as before, then for each generic linear form $l,$ there exists $\epsilon_0$ such that for any $\epsilon$ with $0<\epsilon<\epsilon_0$ and $\delta\neq0$ sufficiently small, the Euler obstruction of $(X,0)$ is equal to 

$$Eu_X(0)=\sum^{q}_{i=1}\chi(V_i\cap B_{\epsilon}\cap l^{-1}(\delta)).Eu_{X}(V_i),$$

\noindent where $\chi$ is the Euler characteristic, $Eu_{X}(V_i)$ is the Euler obstruction of $X$ at a point of $V_i, \ i=1,\ldots,q$ and $0<|\delta|\ll\epsilon\ll1.$
\end{theorem} 


Let us give the definition of another invariant introduced by Brasselet, Massey, Parameswaran and Seade in \cite{BMPS}. Let $f:X\rightarrow\mathbb{C}$ be a holomorphic function with isolated singularity at the origin given by the restriction of a holomorphic function $F:U\rightarrow\mathbb{C}$ and denote by $\overline{\nabla}F(x)$ the conjugate of the gradient vector field of $F$ in $x\in U,$ $$\overline{\nabla}F(x):=\left(\overline{\frac{\partial F}{\partial x_1}},\ldots, \overline{\frac{\partial F}{\partial x_n}}\right).$$

Since $f$ has an isolated singularity at the origin, for all $x\in X\setminus\{0\},$ the projection $\hat{\zeta}_i(x)$ of $\overline{\nabla}F(x)$ over $T_x(V_i(x))$ is nonzero, where $V_i(x)$ is a stratum containing $x.$ Using this projection, the authors constructed, in \cite{BMPS}, a stratified vector field over $X,$ denoted by $\overline{\nabla}f(x).$ Let $\tilde{\zeta}$ be the lifting of $\overline{\nabla}f(x)$ as a section of the Nash bundle $\tilde{T}$ over $\tilde{X}$, without singularity over $\nu^{-1}(X\cap S_{\epsilon}).$

Let $\mathcal{O}(\tilde{\zeta})\in\mathbb{H}^{2n}(\nu^{-1}(X\cap B_{\epsilon}),\nu^{-1}(X\cap S_{\epsilon}))$ be the obstruction cocycle for extending $\tilde{\zeta}$ as a non zero section of $\tilde{T}$ inside $\nu^{-1}(X\cap B_{\epsilon}).$

\begin{definition}
The \textbf{local Euler obstruction of the function} $f, Eu_{f,X}(0)$ is the evaluation of $\mathcal{O}(\tilde{\zeta})$ on the fundamental class $[\nu^{-1}(X\cap B_{\epsilon}),\nu^{-1}(X\cap S_{\epsilon})].$
\end{definition}

The next theorem compares the Euler obstruction of a space $X$ with the Euler obstruction of function defined over $X.$

\begin{theorem}\label{Euler obstruction of a function formula}(Theorem 3.1 of \cite{BMPS})
Let $(X,0)$ and $\mathcal{V}$ be given as before and let \linebreak$f:(X,0)\rightarrow(\mathbb{C},0)$ be a function with an isolated singularity at $0.$ For $0<|\delta|\ll\epsilon\ll1,$ we have
 $$Eu_{f,X}(0)=Eu_X(0)-\sum_{i=1}^{q}\chi(V_i\cap B_{\epsilon}\cap f^{-1}(\delta)).Eu_X(V_i).$$
\end{theorem}



Let us now see a definition we will need to define a generic point of a function-germ. Let $\mathcal{V}=\{V_{\lambda}\}$ be a stratification of a reduced complex analytic space $X.$

\begin{definition}
Let $p$ be a point in a stratum $V_{\beta}$ of $\mathcal{V}.$ A \textbf{degenerate tangent plane of $\mathcal{V}$ at $p$} is an element $T$ of some Grassmanian manifold such that $T=\displaystyle\lim_{p_i\rightarrow p}T_{p_i}V_{\alpha},$ where $p_i\in V_{\alpha}$, $V_{\alpha}\neq V_{\beta}.$
\end{definition}

\begin{definition}
Let $(X,0)\subset(U,0)$ be a germ of complex analytic space in $\mathbb{C}^n$ equipped with a Whitney stratification and let $f:(X,0)\rightarrow(\mathbb{C},0)$ be an analytic function, given by the restriction of an analytic function $F:(U,0)\rightarrow(\mathbb{C},0).$ Then $0$ is said to be a \textbf{generic point}\index{holomorphic function germ!generic point of} of $f$ if the hyperplane $Ker(d_0F)$ is transverse in $\mathbb{C}^n$ to all degenerate tangent planes of the Whitney stratification at $0.$ 
\end{definition}

Now, let us see the definition of a Morsification of a function. 

\begin{definition}
Let $\mathcal{W}=\{W_0,W_1,\ldots,W_q\},$ with $0\in W_0,$ a Whitney stratification of the complex analytic space $X.$ A function $f:(X,0)\rightarrow(\mathbb{C},0)$ is said to be \textbf{Morse stratified} if $\dim W_0\geq1, f|_{W_0}: W_0\rightarrow\mathbb{C}$ has a Morse point at $0$ and $0$ is a generic point of $f$ with respect to $W_{i},$ for all $ i\neq0.$
\end{definition}

A \textbf{stratified Morsification}\index{holomorphic function germ!stratified Morsification of} of a germ of analytic function $f:(X,0)\rightarrow(\mathbb{C},0)$ is a deformation $\tilde{f}$ of $f$ such that $\tilde{f}$ is Morse stratified.

In \cite{STV}, Seade, Tib\u{a}r and Verjovsky proved that the Euler obstruction of a function $f$ is also related to the number of Morse critical points of a stratified Morsification of $f.$

\begin{proposition}(Proposition 2.3 of \cite{STV})\label{Eu_f and Morse points}
Let $f:(X,0)\rightarrow(\mathbb{C},0)$ be a germ of analytic function with isolated singularity at the origin. Then, \begin{center}
$Eu_{f,X}(0)=(-1)^dn_{reg},$
\end{center}
where $n_{reg}$ is the number of Morse points in $X_{reg}$ in a stratified Morsification of $f.$
\end{proposition}

\section{Brasselet number}

\hspace{0,5cm} In this section, we present definitions and results needed in the development of the results of this work. The main reference for this section is \cite{Ms1}.

Let $X$ be a reduced complex analytic space (not necessarily equidimensional) of dimension $d$ in an open set $U\subseteq\mathbb{C}^n$ and let $f:(X,0)\rightarrow(\mathbb{C},0)$ be an analytic map. We write $V(f)=f^{-1}(0).$ 

\begin{definition}\label{good stratification}
A \textbf{good stratification of $X$ relative to $f$} is a stratification $\mathcal{V}$ of $X$ which is adapted to $V(f)$ such that $\{V_{\lambda}\in\mathcal{V},V_{\lambda}\nsubseteq V(f)\}$ is a Whitney stratification of $X\setminus V(f)$ and such that for any pair $(V_{\lambda},V_{\gamma})$ such that $V_{\lambda}\nsubseteq V(f)$ and $V_{\gamma}\subseteq V(f),$ the $(a_f)$-Thom condition is satisfied, that is, if $p\in V_{\gamma}$ and $p_i\in V_{\lambda}$ are such that $p_i\rightarrow p$ and $T_{p_i} V(f|_{V_{\lambda}}-f|_{V_{\lambda}}(p_i))$ converges to some $\mathcal{T},$ then $T_p V_{\gamma}\subseteq\mathcal{T}.$
\end{definition}

If $f:X\rightarrow\mathbb{C}$ has a stratified isolated critical point and $\mathcal{V}$ is a Whitney stratification of $X,$ then \begin{equation}\label{induced stratification}
    \{V_{\lambda}\setminus X^f, V_{\lambda}\cap X^f\setminus\{0\},\{0\}, V_{\lambda}\in\mathcal{V}\}
\end{equation}

\noindent is a good stratification of $X$ relative to $f,$ called the good stratification induced by $f.$

Let $\mathcal{V}$ be a good stratification of $X$ relative to $f.$

\begin{definition}
The \textbf{critical locus of $f$ relative to $\mathcal{V}$}, $\Sigma_{\mathcal{V}}f,$ is given by the union \begin{center}$\Sigma_{\mathcal{V}}f=\displaystyle\bigcup_{V_{\lambda}\in\mathcal{V}}\Sigma(f|_{V_{\lambda}}).$\end{center}
\end{definition}

\begin{definition}
If $\mathcal{V}=\{V_{\lambda}\}$ is a stratification of $X,$ the \textbf{relative polar variety of $f$ and $g$ with respect to $\mathcal{V}$}, denoted by $\Gamma_{f,g}(\mathcal{V}),$ is the the union $\cup_{\lambda}\Gamma_{f,g}(V_{\lambda}),$ where $\Gamma_{f,g}(V_{\lambda})$ denotes the closure in $X$ of the critical locus of $(f,g)|_{V_{\lambda}\setminus X^f},$ where $X^f=X\cap\{f=0\}.$  
\end{definition}

\begin{definition}
If $\mathcal{V}=\{V_{\lambda}\}$ is a stratification of $X,$ the \textbf{symmetric relative polar variety of $f$ and $g$ with respect to $\mathcal{V}$}, $\tilde{\Gamma}_{f,g}(\mathcal{V}),$ is the union $\cup_{\lambda}\tilde{\Gamma}_{f,g}(V_{\lambda}),$ where $\Gamma_{f,g}(V_{\lambda})$ denotes the closure in $X$ of the critical locus of $(f,g)|_{V_{\lambda}\setminus (X^f\cup X^g)},$  $X^f=X\cap \{f=0\}$ and $X^g=X\cap \{g=0\}. $ 
\end{definition}

\begin{definition}\label{definition prepolar}
Let $\mathcal{V}$ be a good stratification of $X$ relative to a function
$f:(X,0)\rightarrow(\mathbb{C},0).$ A function $g :(X, 0)\rightarrow(\mathbb{C},0)$ is \textbf{prepolar with respect to $\mathcal{V}$ at the origin} if the origin is a stratified isolated critical point, that is, $0$ is an isolated point of $\Sigma_{\mathcal{V}}g.$
\end{definition} 




\begin{definition}\label{definition tractable}
A function $g :(X, 0)\rightarrow(\mathbb{C},0)$ is \textbf{tractable at the origin with respect to a good stratification $\mathcal{V}$ of $X$ relative to $f :(X, 0)\rightarrow(\mathbb{C},0)$} if $dim_0 \ \tilde{\Gamma}^1_{f,g}(\mathcal{V})\leq1$ and, for all strata $V_{\alpha}\subseteq X^f$,
$g|_{V_{\alpha}}$ has no critical point in a neighbourhood of the origin except perhaps at the origin itself.

\end{definition}

Let us now see the definition of decent analytic function-germs. Let $\mathcal{V}=\{V_{\lambda}\}$ be a stratification of a reduced complex analytic space $X.$

\begin{definition}
Let $g:(X,0)\rightarrow(\mathbb{C},0)$ be a function-germ. For any analytic stratification $\mathcal{V}$ of $X$, $0<|\delta|\ll1$ and any function $f: (X, 0)\rightarrow (\mathbb{C},0),$ $g$ is \textbf{decent with respect to $\mathcal{V}$ relative to $f$} if there exists a neighborhood $\Omega$ of $0$ such that $g:\Omega\cap X\cap f^{-1}(\delta)\setminus X^g\rightarrow\mathbb{C}$ has only generic points.
\end{definition}

\begin{proposition}(Proposition 1.14 of \cite{Ms1})\label{Proposition 1.14}
Let $\mathcal{V}$ be a good stratification of $X$ relative to $f$ at the origin. Then, for a generic choice of linear form, $l,$ $l$ is decent to $\mathcal{V}$ relative to $f$ and, moreover, $f$ is decent with respect to $\mathcal{V}$ relative to $l.$
\end{proposition}

Another concept useful for this work is the notion of constructible functions. Consider a Whitney stratification $\mathcal{W}=\{W_1,\ldots, W_q\}$ of $X$ such that each stratum $W_i$ is connected.

\begin{definition}
A constructible function with respect to the stratification $\mathcal{W}$ of $X$ is a function $\beta:X\rightarrow\mathbb{Z}$ which is constant on each stratum $W_i,$ that is, there exist integers $t_1,\ldots,t_q,$ such that $\beta=\sum_{i=1}^{q}t_i.1_{W_i},$ where $1_{W_i}$ is the characteristic function of $W_i.$
\end{definition}

\begin{definition}
The Euler characteristic $\chi(X,\beta)$ of a constructible function $\beta:X\rightarrow\mathbb{Z}$ with respect to the stratification $\mathcal{W}$ of $X,$ given by $\beta=\sum_{i=1}^{q}t_i.1_{W_i}$, is defined by $\chi(X,\beta)=\sum_{i=1}^{q}t_i.\chi(W_i).$
\end{definition}



Before we state Dutertre and Grulha results, we need to introduce some definitions about normal Morse data. We cite as main references \cite{goresky1988stratified} and \cite{schurmann2010index}. The first concept we present is the complex link, an object analogous to the Milnor fibre, important in the study of complex stratified Morse theory. 

Let $V$ be a stratum of the stratification $\mathcal{V}$ of $X$ and let $x$ be a point of $V.$ Let \linebreak$g:(\mathbb{C}^n,0)\rightarrow(\mathbb{C},0)$ be an analytic complex function-germ such that the differential form $Dg(x)$ does not vanish on a degenerate tangent plane of $\mathcal{V}$ at $x.$ Let $N$ be a normal slice to $V$ at $x,$ that is, $N$ is a closed complex submanifold of $\mathbb{C}^n$ which is transversal to $V$ at $x$ and $N\cap V=\{x\}.$

\begin{definition}
Let $B_{\epsilon}$ be the closed ball of radius $\epsilon$ centered at $x.$ The \textbf{complex link} $l_V$ of $V$ is defined by $l_V=X\cap N\cap B_{\epsilon}\cap\{g=\delta\},$ where $0<|\delta|\ll\epsilon\ll1.$

The \textbf{normal Morse datum} $NMD(V)$ of $V$ is the pair of spaces \begin{eqnarray*}
NMD(V)=(X\cap N\cap B_{\epsilon},X\cap N\cap B_{\epsilon}\cap\{g=\delta\}).
\end{eqnarray*}
 \end{definition}
 
In Part II, section 2.3 of \cite{goresky1988stratified}, the authors explained why this two notions are independent of all choices made. 

\begin{definition}
 Let $\beta:X\rightarrow\mathbb{Z}$ be a constructible function with respect to the stratification $\mathcal{V}.$ Its normal Morse index $\eta(V,\beta)$ along $V$ is defined by 
 \begin{eqnarray*}
 \eta(V,\beta)=\chi(NMD(V),\beta)=\chi(X\cap N\cap B_{\epsilon},\beta)-\chi(l_V,\beta).
 \end{eqnarray*}
\end{definition}
 
In the case where the constructible function is the local Euler obstruction, the following identities are valid (\cite{schurmann2010index}, page 34):
\begin{center}
    $\eta(V^{\prime},Eu_{\overline{V}})=1,$ if $V^{\prime}=V$ and $\eta(V^{\prime},Eu_{\overline{V}})=0,$ if $V^{\prime}\neq V.$
\end{center}

We present now the definition of the Brasselet number and the main theorems of \cite{DG}, used as inspiration for this work.

Let $f: (X,0)\rightarrow(\mathbb{C},0)$ be a complex analytic function germ and let $\mathcal{V}$ be a good stratification of $X$ relative to $f.$ We denote by $V_1,\ldots, V_q$ the strata of $\mathcal{V}$ that are not contained in $\{f=0\}$ and we assume that $V_1,\ldots, V_{q-1}$ are connected and that $V_{q}=\linebreak X_{reg}\setminus \{f=0\}.$ Note that $V_q$ could be not connected.  
 
\begin{definition}
Suppose that $X$ is equidimensional. Let $\mathcal{V}$ be a good stratification of $X$ relative to $f.$ The \textbf{Brasselet number} of $f$ at the origin, $B_{f,X}(0),$ is defined by \begin{center}
$B_{f,X}(0)=\sum_{i=1}^{q}\chi(V_i\cap f^{-1}(\delta)\cap B_{\epsilon})Eu_X(V_i),$
\end{center}
where $0<|\delta|\ll\epsilon\ll1.$
\end{definition} 

\noindent\textbf{Remark:} If $V_q^i$ is a connected component of $V_{q},$ $Eu_X(V_q^i)=1.$

Notice that if $f$ has a stratified isolated singularity at the origin, then \linebreak$B_{f,X}(0)=Eu_{X}(0)-Eu_{f,X}(0)$ (see Theorem \ref{Euler obstruction of a function formula}).

In \cite{DG}, Dutertre and Grulha proved interesting formulas describing the topological relation between the Brasselet number and a number of certain critical points of a special type of deformation of functions. Let us now present some of these results. Fist we need the definition of a special type of Morsification, introduced by Dutertre and Grulha.

\begin{definition}
A \textbf{partial Morsification} of $g:f^{-1}(\delta)\cap X\cap B_{\epsilon}\rightarrow\mathbb{C}$ is a function $\tilde{g}: f^{-1}(\delta)\cap X\cap B_{\epsilon}\rightarrow\mathbb{C}$ (not necessarily holomorphic) which is a local Morsification of all isolated critical points of $g$ in $f^{-1}(\delta)\cap X\cap \{g\neq 0\}\cap B_{\epsilon}$ and which coincides with $g$ outside a small neighborhood of these critical points.
\end{definition} 

Let $g : (X, 0) \rightarrow (\mathbb{C},0)$ be a complex analytic function which is tractable at the origin with respect to $\mathcal{V}$ relative to $f .$ Then $\tilde{\Gamma}_{f,g}$ is a complex analytic curve and for $0<|\delta|\ll1$ the critical points of $g|_{f^{-1}(\delta)\cap X}$ in $B_{\epsilon}$ lying outside $\{g=0\}$ are isolated.
Let $\tilde{g}$ be a partial Morsification of $g:f^{-1}(\delta)\cap X\cap B_{\epsilon}\rightarrow\mathbb{C}$ and, for each \linebreak$i\in\{1,\ldots, q\},$ let $n_i$ be the number of stratified Morse critical points of $\tilde{g}$ appearing on \linebreak$V_i\cap f^{-1}(\delta)\cap \{g\neq 0\}\cap B_{\epsilon}.$ 
 
\begin{theorem}\label{4.2 DG}(Theorem 4.2 of \cite{DG})
Let $\beta: X\rightarrow\mathbb{Z}$ be a constructible function with respect to the stratification $\mathcal{V}$. Suppose that $g:(X,0)\rightarrow(\mathbb{C},0)$ is a complex analytic function tractable at the origin with respect to $\mathcal{V}$ relative to $f.$ For $0<|\delta|\ll\epsilon\ll 1,$ we have\begin{center}
 $\chi(X\cap f^{-1}(\delta)\cap B_{\epsilon},\beta)-\chi(X\cap g^{-1}(0)\cap f^{-1}(\delta)\cap B_{\epsilon},\beta)=\sum_{i=1}^{q}(-1)^{d_i-1}n_i\eta(V_i,\beta).$
 \end{center}
\end{theorem}

In the case that $\beta=Eu_X$, the last theorem implies the following.

\begin{corollary}\label{Corollary 4.3 of DG}(Corollary 4.3 of \cite{DG})
Suppose that $X$ is equidimensional and that $g$ is tractable at the origin with respect to $\mathcal{V}$ relative to $f.$ For $0<|\delta|\ll\epsilon\ll 1,$ we have\begin{center}
 $\chi(X\cap f^{-1}(\delta)\cap B_{\epsilon},Eu_X)-\chi(X\cap g^{-1}(0)\cap f^{-1}(\delta)\cap B_{\epsilon},Eu_X)=(-1)^{d-1}n_q.$
 \end{center}
\end{corollary}

If one supposes, in addition, that $g$ is prepolar, a consequence of this result is a Lê-Greuel type formula for the Brasselet number.

\begin{theorem}\label{4.4 DG}(Theorem 4.4 of \cite{DG})
Suppose that $X$ is equidimensional and that $g$ is prepolar with respect to $\mathcal{V}$ at the origin. For $0<|\delta|\ll\epsilon\ll 1,$ we have\begin{center}
 $B_{f,X}(0)-B_{f,X^g}(0)=(-1)^{d-1}n_q,$
 \end{center}
 where $n_q$ is the number of stratified Morse critical points on the top stratum $V_q\cap f^{-1}(\delta)\cap B_{\epsilon}$ appearing in a Morsification of $g:X\cap f^{-1}(\delta)\cap B_{\epsilon}\rightarrow \mathbb{C}.$
\end{theorem}


Suppose that $X$ is equipped with a Whitney stratification $\mathcal{V}=\{V_0,V_1,\ldots,V_q\}$ with $V_0=\{0\},$ and $f,g:X\rightarrow\mathbb{C}$ have an isolated stratified singularity at the origin with respect to this stratification. We give now some results proved by Dutertre and Grulha in Section 6 of \cite{DG} in this setting.

\begin{proposition}
Suppose that $g$ (resp. $f$) is prepolar with respect to the good stratification induced by $f$ (resp. $g$) at the origin. Let $\beta:X\rightarrow\mathbb{Z}$ be a constructible function with respect to the Whitney stratification $\mathcal{V}.$ For $0<|\delta|\ll\epsilon\ll1,$ \begin{center}
$\chi(X^f\cap g^{-1}(\delta)\cap B_{\epsilon},\beta)=\chi(X^g\cap f^{-1}(\delta)\cap B_{\epsilon},\beta).$
\end{center}
\end{proposition}

A corollary of this proposition is the following result.

\begin{corollary}\label{corollary 6.3 DG}
Suppose that $X$ is equidimensional and that $g$ (resp. $f$) is prepolar with respect to the good stratification induced by $f$ (resp. $g$) at the origin. Then \begin{center}$B_{f,X^g}(0)=B_{g,X^f}(0).$\end{center}
\end{corollary}

In \cite{DG}, the authors also related the topology of the  generalized Minor fibres of $f$ and $g$ and some number of Morse points. 
\begin{theorem}\label{6.4 of DG}
Suppose that $g$ (resp. $f$) is prepolar with respect to the good stratification induced by $f$ (resp. $g$) at the origin. Let $\beta:X\rightarrow\mathbb{Z}$ be a constructible function with respect to the Whitney stratification $\mathcal{V}.$ For $0<|\delta|\ll\epsilon\ll1,$ \begin{center}
$\chi(X\cap f^{-1}(\delta)\cap B_{\epsilon},\beta)-\chi(X\cap g^{-1}(\delta)\cap B_{\epsilon},\beta)=\sum_{i=1}^{q}(-1)^{d_i-1}(n_i-m_i)\eta(V_i,\beta),$
\end{center}
\noindent where $n_i$ (resp. $m_i$) is the number of stratified Morse critical points on the stratum \linebreak$V_i\cap f^{-1}(\delta)\cap B_{\epsilon}$ (resp. $V_i\cap g^{-1}(\delta)\cap B_{\epsilon}$) appearing in a Morsification of \linebreak$g:X\cap f^{-1}(\delta)\cap B_{\epsilon}\rightarrow\mathbb{C}$ (resp. $f:X\cap g^{-1}(\delta)\cap B_{\epsilon}\rightarrow\mathbb{C}$).
\end{theorem}

In the case where $\beta=Eu_X,$ the last theorem implies the following result.

\begin{corollary}\label{corollary 6.5 of DG}
Suppose that $X$ is equidimensional and that $g$ (resp. $f$) is prepolar with respect to the good stratification induced by $f$ (resp. $g$) at the origin. Then \begin{center}
$B_{f,X}(0)-B_{g,X}(0)=(-1)^{d-1}(n_q-m_q),$
\end{center}
\noindent where $n_q$ (resp. $m_q$) is the number of stratified Morse critical points on the top stratum \linebreak$V_q\cap f^{-1}(\delta)\cap B_{\epsilon}$ (resp. $V_q\cap g^{-1}(\delta)\cap B_{\epsilon}$) appearing in a Morsification of \linebreak$g:X\cap f^{-1}(\delta)\cap B_{\epsilon}\rightarrow\mathbb{C}$ (resp. $f:X\cap g^{-1}(\delta)\cap B_{\epsilon}\rightarrow\mathbb{C}$).
\end{corollary}

Applying Corollary \ref{corollary 6.3 DG} to the case where the function $g$ is a generic linear form, one obtains the following result.
\begin{corollary}\label{6.6 of DG}
Suppose that $X$ is equidimensional. Let $H$ be a generic hyperplane. Then \begin{center}
$Eu_{X^f}(0)=B_{f,X\cap H}(0).$
\end{center}
\end{corollary}

The Brasselet number $B_{f,X\cap H}(0)$ can also be compared to $B_{f,X}(0)$ using the dimension $d$ of $(X,0)$ and the generic polar curves. Let $\Gamma^0_{f|_X})$ be the general relative polar curve of the morphism $f:X\rightarrow\mathbb{C}$ (see \cite{LT} for the definition of the general relative polar curve). Let us denote by $I_0(X^f,\Gamma^0_{f|_X})$ the cardinality of $\Gamma^0_{f|_X}\cap f^{-1}(\delta  ), 0\ll|\delta|\ll1,$ that is, the intersection multiplicity of $\Gamma^0_{f|_X}$ at $0$ in $X^f$ (see \cite{loeser1984formules}).

\begin{corollary}\label{corollary 5.2 of DG}
Suppose that $X$ is equidimensional. Then \begin{center}
    $B_{f,X}(0)-B_{f,X\cap H}(0)=(-1)^{d-1}I_0(X^f,\Gamma^0_{f|_X}),$
\end{center}
\noindent where $H$ is a generic hyperplane.
\end{corollary}

If the symmetric relative polar variety $\tilde{\Gamma}_{f,g}(\mathcal{V})$ has dimension one, for each $V_i\in\mathcal{V},$ one associates a \textit{multiplicity} $\mu^f(\tilde{\Gamma}_{f,g}(V_i))$ to each $\tilde{\Gamma}_{f,g}(V_i)$ in the following way: if the stratum $V_i$ is one-dimensional, $\mu^f(\tilde{\Gamma}_{f,g}(V_i)):=1.$ If $V_i$ is not one-dimensional, let $\nu$ be a component of $\tilde{\Gamma}_{f,g}(V_i)$ and $p$ be a point of $\nu\setminus\{0\}$ close to the origin. The mapping $g:V_i\cap\{f=f(p)\}\rightarrow\mathbb{C}$ has an isolated singularity at $p$ and let $\mu^{\nu}$ be the Milnor number of this singularity. Then $\mu^f(\tilde{\Gamma}_{f,g}(V_i))$ is the sum of Milnor numbers $\mu^{\nu}$ over all components $\nu.$ (See page 974 of \cite{Ms1})

At last, a consequence of Corollary \ref{corollary 6.5 of DG} in the case where $g$ is a generic linear form is the following.

\begin{corollary}\label{6.7 of DG}
Suppose that $X$ is equidimensional. Let $l$ be a generic linear form. Then \begin{center}
$\mu^f(\Gamma_{f,l}(V_q))-\mu^l(\Gamma_{f,l}(V_q))=(-1)^dEu_{f,X}(0).$
\end{center}

\end{corollary}


 \section{Some results for functions with arbitrary singularities}

\hspace{0,5cm} Let $(X,0)$ be a reduced equidimensional analytic germ of dimension $d$ in an open set $U\subset\mathbb{C}^n$ and $f,g:(X,0)\rightarrow(\mathbb{C},0)$ be two germs of functions. Let $\mathcal{V}$ be a good stratification of $X$ relative to $f$ and suppose that the critical locus of $g,$ $\Sigma_{\mathcal{V}} g,$ is one-dimensional and that $\Sigma_{\mathcal{V}} g\cap \{f=0\}=\{0\}.$

Let $V_1,\ldots, V_q$ be the strata of $\mathcal{V}$ not contained in $\{f=0\}.$ Suppose that $\{0\}$ is a stratum of $\{f=0\}$, that for each $i\in\{1,\ldots,q-1\},$ $V_i$ is connected, $V_q$ is equal to $X_{reg}\setminus\{f=0\}$ and that $d_i=\dim V_i.$ In this case, we can construct a good stratification of $X$ relative to $f$ that gives us also a good stratification of $X\cap\{g=0\}$ relative to $f|_{X\cap\{g=0\}}.$ We start this section with the construction of this stratification.

\begin{lemma}\textbf{(First stratification lemma)}\label{first stratification lemma}
Let $\mathcal{V}$ be a good stratification of $X$ relative to $f$ and $\mathcal{V}^f$ be the collection of strata of $\mathcal{V}$ contained in $\{f=0\}$ (including the stratum $\{0\}$). Then, the refinement 
\begin{eqnarray}
\mathcal{V}^{\prime}=\Big\{ V_{i}\setminus \Sigma_{\mathcal{V}} g, V_{i}\cap \Sigma_{\mathcal{V}} g, i\in\{1,\ldots,q\}\Big\} \cup\mathcal{V}^f\label{good stratification}
\end{eqnarray} is a good stratification of $X$ relative to $f$ and  $\mathcal{V}^{\prime\{g=0\}}$ is a good stratification of $X\cap\{g=0\}$ relative to $f|_{X\cap\{g=0\}},$ where \begin{center}$\mathcal{V}^{\prime\{g=0\}}=\Big\{V_{i}\cap\{g=0\} \setminus \Sigma_{\mathcal{V}} g, V_{i}\cap \Sigma_{\mathcal{V}} g, i\in\{1,\ldots,q\}\Big\}\cup\mathcal{V}^f\cap\{g=0\},$\end{center}
\noindent and $\mathcal{V}^f\cap\{g=0\}$ denotes the collection of strata of type $V^f\cap\{g=0\},$ with $V^f\in\mathcal{V}^f.$
Moreover, if $g$ is tractable at the origin with respect to $\mathcal{V}$ relative to $f,$ then $g$ is tractable at the origin with respect to $\mathcal{V}^{\prime}$ relative to $f.$
\end{lemma}

\noindent\textbf{Proof.} Since $\Sigma_{\mathcal{V}} g\cap \{f=0\}=\{0\}$ and $V_1,\ldots, V_q$ are the strata of $\mathcal{V}$ not contained in $\{f=0\},$ we can write $\Sigma_{\mathcal{V}} g=\{0\}\cup(V_1\cap\Sigma_{\mathcal{V}} g)\cup\ldots\cup (V_q\cap \Sigma_{\mathcal{V}} g).$ Let us show that the refinement of $\mathcal{V},$
\begin{equation*}
\mathcal{V}^{\prime}=\Big\{ V_{i}\setminus \Sigma_{\mathcal{V}} g, V_{i}\cap \Sigma_{\mathcal{V}} g, i\in\{1,\ldots,q\}\Big\}\cup\mathcal{V}^f,
\end{equation*}  
is a good stratification of $X$ relative to $f.$ Since the collection of strata contained in $\{f=0\}$ was not refined, $\{f=0\}$ is union of strata of $\mathcal{V}^{\prime}.$ Now we will show that \begin{center}$\Big\{V_{\alpha}\in\mathcal{V}^{\prime}; V_{\alpha}\nsubseteq\{f=0\}\Big\}=\Big\{V_{i}\setminus \Sigma_{\mathcal{V}} g, V_{i}\cap \Sigma_{\mathcal{V}} g, i\in\{1,\ldots,q\}\Big\}$\end{center}
is a Whitney stratification of $X\setminus\{f=0\}.$ We can refine this stratification to obtain a Whitney stratification. But since $\mathcal{V}$ is a good stratification of $X$ relative to $f,$ \linebreak$\{V_{\alpha}\in\mathcal{V}; V_{\alpha}\nsubseteq\{f=0\}\}$ is a Whitney stratification of $X\setminus\{f=0\}.$ Since $\Sigma_{\mathcal{V}} g$ is closed, $V_{\alpha}\setminus\Sigma_{\mathcal{V}} g$ is an open subset of $V_{\alpha}$ and then the Whitney's condition $(b)$ is verified over the strata of type $V_{\alpha}\setminus\Sigma_{\mathcal{V}} g.$ So, the refinement should be done only over the stratum of type $\Sigma_{\mathcal{V}} g\cap V_i.$ Since $\Sigma_{\mathcal{V}} g$ is one-dimensional, a refinement of $\Sigma_{\mathcal{V}} g\cap V_i$ would be done by taking off a finite number of points. So, in a sufficiently small neighborhood of the origin, the Whitney's condition $(b)$ is verified over $\{V_{i}\setminus \Sigma_{\mathcal{V}} g, V_{i}\cap \Sigma_{\mathcal{V}} g, i\in\{1,\ldots,q\}\}.$

At last, let us verify the Thom condition. Let $p$ be a point in $V_{\beta}\subseteq \{f=0\}$ and $(p_k)$ be a sequence of points in $V_{\alpha}\nsubseteq \{f=0\}.$ Suppose that $\displaystyle\lim_{k\rightarrow\infty}p_k= p$ and that\linebreak $\displaystyle\lim_{k\rightarrow\infty}T_{p_k}V(f|_{V_{\alpha}}-f|_{V_{\alpha}}(p_k))=T.$ We must show that $T_p V_{\beta}\subseteq T.$ 

Let $V_{\alpha}=V_i\cap\Sigma_{\mathcal{V}} g\,$ for some $i\in\{1,\ldots, q\}.$ 
Since $(p_k)$ is a sequence of points in $V_i\cap\Sigma_{\mathcal{V}} g$ and $\Sigma_{\mathcal{V}} g$ is one-dimensional, $(p_k)$ must converge to the origin, that is, $p=0$ and $V_{\beta}=\{0\}.$ So, $\{0\}=T_p V_{\beta}\subseteq T.$

Let us now to verify the Thom condition for $V_{\alpha}=V_i\setminus \Sigma_{\mathcal{V}} g, i\in\{1,\ldots,q\}.$ We have \begin{center}
$T_{p_k}V(f|_{{V_i\setminus \Sigma_{\mathcal{V}} g}}-f|_{{V_i\setminus \Sigma_{\mathcal{V}} g}}(p_k))=T_{p_k}V(f|_{V_i}-f|_{V_i}(p_k)),$
\end{center} 
and since $\mathcal{V}$ is a good stratification of $X$ relative to $f,$ the Thom condition is verified for $(V_i, V_{\beta}).$ Hence, $T_p V_{\beta}\subseteq \displaystyle \lim_{k\rightarrow\infty}T_{p_k}V(f|_{V_i}-f|_{V_i}(p_k))=\lim_{k\rightarrow\infty}T_{p_k}V(f|_{{V_i\setminus \Sigma_{\mathcal{V}} g}}-f|_{{V_i\setminus \Sigma_{\mathcal{V}} g}}(p_k))=T.$ Therefore, $\mathcal{V}^{\prime}$ is a good stratification of $X$ relative to $f.$

Let us now show that $\mathcal{V}^{\prime\{g=0\}}$ is a good stratification of $X\cap\{g=0\}$ relative to $f|_{X\cap\{g=0\}},$ \begin{center}$\mathcal{V}^{\prime\{g=0\}}=\Big\{V_{i}\cap\{g=0\} \setminus \Sigma_{\mathcal{V}} g, V_{i}\cap \Sigma_{\mathcal{V}} g, i\in\{1,\ldots,q\}\Big\}\cup\mathcal{V}^f\cap\{g=0\}.$\end{center}
Since $ \Sigma_{\mathcal{V}} g\cap\{f=0\}=\{0\}, \{g=0\}$ intersects each stratum of $\mathcal{V}^f$ transversely. Therefore, for each $\tilde{V_i}\in \mathcal{V}^f,$ $\{g=0\}\cap \tilde{V_i}$ is a complex analytic submanifold of $\tilde{V_i}.$ Hence, $\{f=0\}\cap\{g=0\}$ is union of strata contained in $\mathcal{V}^f\cap\{g=0\}.$ Now, we will verify that \begin{center}$\{V_{i}\cap\{g=0\} \setminus \Sigma_{\mathcal{V}} g, V_{i}\cap \Sigma_{\mathcal{V}} g, i\in\{1,\ldots,q\}\}$\end{center} stratification of $X\cap\{g=0\}\setminus\{f=0\}.$ Consider a pair of strata of type \linebreak $(V_i\cap\Sigma_{\mathcal{V}}g,V_j\cap\Sigma_{\mathcal{V}}g).$ If necessary, we can refine these strata to guarantee the Whitney's condition $(b)$. Since $\Sigma_{\mathcal{V}}g$ has dimension one, this refinement would be given by taking off a finite number of points. Therefore, in a sufficiently small neighborhood of the origin, Whitney's condition $(b)$ is verified for this type of stratum. Now, let us verify this condition for pairs of strata of the type $(V_{i}\cap\{g=0\} \setminus \Sigma_{\mathcal{V}} g,V_{j}\cap\{g=0\} \setminus \Sigma_{\mathcal{V}} g)$ and $(V_{i}\cap\{g=0\} \setminus \Sigma_{\mathcal{V}} g,V_{j}\cap\Sigma_{\mathcal{V}} g).$ 
\begin{enumerate}
    \item Let us show that $(V_{i}\cap\{g=0\} \setminus \Sigma_{\mathcal{V}} g,V_{j}\cap\{g=0\} \setminus \Sigma_{\mathcal{V}} g)$ is Whitney regular. Since these strata contain no critical points of $g, V_{i}\cap\{g=0\}$ and $V_{j}\cap\{g=0\}$ are transverse intersections. Therefore, $(V_{i}\cap\{g=0\} \setminus \Sigma_{\mathcal{V}} g,V_{j}\cap\{g=0\} \setminus \Sigma_{\mathcal{V}} g)$ is Whitney regular, since $(V_i,V_j)$ is Whitney regular. (See \cite{orro2010regularity}).

    \item Let us show that $(V_{i}\cap\{g=0\} \setminus \Sigma_{\mathcal{V}} g,V_{j}\cap\Sigma_{\mathcal{V}} g)$ is Whitney regular. The intersection $V_{i}\cap\{g=0\}$ is transverse, since it contains no critical points of $g.$ Whitney's condition $(b)$ could fail over $V_{j}\cap\Sigma_{\mathcal{V}} g,$ but since $\Sigma_{\mathcal{V}}g$ is one-dimensional, we can refine this stratum by taking off a finite number of points and ensure that, in a sufficiently small neighborhood of the origin, $(V_{i}\cap\{g=0\} \setminus \Sigma_{\mathcal{V}} g,V_{j}\cap\Sigma_{\mathcal{V}} g)$ is Whitney regular.
    
\end{enumerate}
Let us now verify the Thom condition over the strata of $\mathcal{V}^{\prime\{g=0\}}.$ Let $V_{\alpha}\nsubseteq \{f=0\}$ and $V_{\beta}\subseteq \{f=0\}$ be strata of $\mathcal{V}^{\prime\{g=0\}},$ $p$ be a point in $V_{\beta}$ and $(p_i)$ be a sequence of points in $V_{\alpha}.$ Suppose that $\displaystyle\lim_{i\rightarrow\infty}(p_i)= p$ and that $\displaystyle\lim_{i\rightarrow\infty}T_{p_i}V(f|_{V_{\alpha}}-f|_{V_{\alpha}}(p_i))=T.$ We must show that $T_p V_{\beta}\subseteq T.$ If $V_{\beta}=\{0\},$ $T_p V_{\beta}=\{0\}$ and then $T_p V_{\beta}\subseteq T.$ Suppose now that $V_{\beta}\neq\{0\}.$ Notice that, since $\Sigma_{\mathcal{V}}g\cap \{f=0\}=\{0\}$ and $\{0\}\neq V_{\beta}\subseteq\{f=0\}, p\in V_{\beta}$ implies that $p\not\in\Sigma_{\mathcal{V}}g.$ As we have seen above, it is sufficient to verify the Thom condition for $V_{\alpha}=V_j\cap \{g=0\}\setminus \Sigma_{\mathcal{V}} g.$ We have, \begin{eqnarray*}T_{p_i}V(f|_{V_{\alpha}}-f|_{V_{\alpha}}(p_i))&=&T_{p_i}V(f|_{V_{j}\cap \{g=0\}\setminus\Sigma_{\mathcal{V}} g}-f|_{V_{j}\cap \{g=0\}\setminus\Sigma_{\mathcal{V}} g}(p_i))\\&=&T_{p_i}V(\tilde{f}|_{V_{j}}-\tilde{f}|_{V_{j}}(p_i))\cap T_{p_i} V(\tilde{g}),\end{eqnarray*}
which implies that \begin{eqnarray*}
\displaystyle\lim_{i\rightarrow\infty}T_{p_i}V(f|_{V_{\alpha}}-f|_{V_{\alpha}}(p_i))&=&\displaystyle\lim_{i\rightarrow\infty} T_{p_i}V(\tilde{f}|_{V_{j}}-\tilde{f}|_{V_{j}}(p_i))\cap T_{p_i}V(\tilde{g})\\&\subseteq& \lim_{i\rightarrow\infty} T_{p_i}V(\tilde{f}|_{V_{j}}-\tilde{f}|_{V_{j}}(p_i))\cap \lim_{i\rightarrow\infty} T_{p_i}V(\tilde{g}) ,
\end{eqnarray*}
where $\tilde{f}$ and $\tilde{g}$ denote analytic extensions of $f$ and $g$ to an open neighborhood of the origin in the ambient space $(U,0)$ of $(X,0).$
Since $p, p_i\not\in\Sigma_{\mathcal{V}} g,$ if $\displaystyle\lim_{i\rightarrow\infty}T_{p_i}V(\tilde{f}|_{V_{j}}-\tilde{f}|_{V_{j}}(p_i))=T_1,$ we have that $T\subseteq T_1\cap T_p V(\tilde{g}).$
On the other hand, if $V_{\beta}=\tilde{V}_{\lambda}\cap\{g=0\},$ with $\tilde{V}_{\lambda}\in\mathcal{V}^f$ and $p\not\in\Sigma_{\mathcal{V}}g,$ \begin{center}$T_pV_{\beta}=T_p(\tilde{V_{\lambda}}\cap \{g=0\})= T_p\tilde{V_{\lambda}}\cap T_pV(\tilde{g}).$\end{center} Since the Thom condition is valid over $\mathcal{V}^{\prime},$ we have that $T_p\tilde{V_{\lambda}}\subseteq T_1.$ Since $g$ is tractable at the origin with respect to $\mathcal{V}$ relative to $f,$ for $p\not\in\Sigma_{\mathcal{V}}g, T_pV(\tilde{g})$ intersects $T_p\tilde{V}_{\lambda}$ transversely. Therefore, $T_pV(\tilde{g})$ intersects $T_1$ transversely. This implies that \begin{eqnarray*}
T=\displaystyle\lim_{i\rightarrow\infty} T_{p_i}V(\tilde{f}|_{V_{j}}-\tilde{f}|_{V_{j}}(p_i))\cap T_{p_i}V(\tilde{g})&=& \lim_{i\rightarrow\infty} T_{p_i}V(\tilde{f}|_{V_{j}}-\tilde{f}|_{V_{j}}(p_i))\cap \lim_{i\rightarrow\infty} T_{p_i}V(\tilde{g})\\&=& T_1\cap T_pV(\tilde{g}).
\end{eqnarray*}

Therefore,\begin{center}
$T_pV_{\beta}=T_p(\tilde{V_{\lambda}}\cap \{g=0\})= T_p\tilde{V_{\lambda}}\cap T_pV(\tilde{g})\subseteq T_1\cap T_pV(\tilde{g})=T.$
\end{center}

Suppose now that $g$ is tractable at the origin with respect to $\mathcal{V}$ relative to $f.$ Let us show that $g$ is tractable at the origin with respect to $\mathcal{V}^{\prime}$ relative to $f.$ For that we should verify that (1): $\dim_0\tilde{\Gamma}_{f,g}(\mathcal{V}^{\prime})\leq1;$ and that (2): $g|_{V_{\alpha}}$ has no singularity in a neighborhood of the origin, except perhaps the origin itself, for $V_{\alpha}\in\mathcal{V}^{\prime}$ contained in $\{f=0\}$. Condition (2) is valid, since we have not refined the strata contained in $\{f=0\}.$ Let us verify condition (1). Since for each $V_i\nsubseteq\{f=0\}, \Sigma g|_{V_i}\subset \{g=0\},$ we have \begin{center}
    $\displaystyle\tilde{\Gamma}_{f,g}(\mathcal{V}^{\prime})=\bigcup_{i=1}^{q}\overline{\Sigma(f,g)|_{(V_i\setminus\Sigma g|_{V_i})\setminus\{f=0\}\cup\{g=0\}}}=\bigcup_{i=1}^{q}\overline{\Sigma(f,g)|_{V_i\setminus\{f=0\}\cup\{g=0\}}}=\tilde{\Gamma}_{f,g}(\mathcal{V}).$
\end{center}
Then $\dim_0 \tilde{\Gamma}_{f,g}(\mathcal{V}^{\prime})=\dim_0 \tilde{\Gamma}_{f,g}(\mathcal{V})\leq1$ and condition (1) is verified. Therefore, $g$ is tractable at the origin with respect to $\mathcal{V}^{\prime}.$
\fim

\vspace{0,3cm}

 By definition, $\Sigma_{\mathcal{V}} g=\bigcup_{\alpha=1}^q\Sigma g|_{V_{\alpha}}\cup\{0\},$ where $V_{\alpha}$ is a stratum not contained in $\{f=0\}.$ Since $\Sigma_{\mathcal{V}} g$ is one-dimensional at the origin, $\Sigma g|_{V_{\alpha}}$ is either one-dimensional or the origin itself, and $\overline{\Sigma g|_{V_{\alpha}}}=\Sigma g|_{V_{\alpha}}\cup\{0\}.$ Since $\overline{\Sigma g|_{V_{\alpha}}}=\overline{\Sigma g|_{\overline{V_{\alpha}}}\setminus Sing(\overline{V_{\alpha}})}$ is an analytic set, it has an irreducible decomposition into one-dimensional subvarieties,  which will be called branches,
\begin{eqnarray*}
\overline{\Sigma g|_{V_{\alpha}}}=\Sigma g|_{V_{\alpha}}\cup\{0\}=b_{\alpha_1}\cup\ldots\cup b_{\alpha_t}.
\end{eqnarray*}

Making this process for each stratum $V_{\alpha},$ we can decompose $\Sigma_{\mathcal{V}}g$ into branches $b_j,$
\begin{eqnarray*}
\Sigma_{\mathcal{V}} g=\bigcup_{\alpha=1}^q\Sigma g|_{V_{\alpha}}\cup\{0\}=b_1\cup\ldots\cup b_r,
\end{eqnarray*}

\noindent where $b_j\subseteq V_{\alpha},$ for some $\alpha\in\{1,\ldots,q\}.$ Notice that a stratum $V_{\alpha}$ can contain no branch and that a stratum  $V_j$ can contain more than one branch, but, the way we described, a branch can not be contained in two different strata. Let $\delta$ be a regular value of $f, 0<|\delta|\ll1,$ and let us write, for each $j\in\{1,\ldots,r\}, f^{-1}(\delta)\cap b_j=\{x_{i_1},\ldots,x_{i_{k(j)}}\}.$ So, in this case, the local degree $m_{f,b_j}$ of $f|_{b_j}$ is $k(j).$ Let $\epsilon$ be sufficiently small such that the local Euler obstruction of $X$ and of $X^g$ are constant on $b_j\cap B_{\epsilon}$. In this case, we denote by $Eu_{X}(b_j)$ (respectively, $Eu_{X^g}(b_j)$) the local Euler obstruction of $X$ (respectively, $X^g$) at a point of $b_j\cap B_{\epsilon}.$

The next theorem calculates, in our setting, the difference $B_{f,X}(0)-B_{f,X^g}(0)$ without the prepolarity of $g$ with respect to the good stratification relative to $f$ at the origin.

We fix the good stratification $\mathcal{V}^{\prime}$ of $X$ relative to $f$ constructed in Lemma \ref{first stratification lemma} given as a refinement of the initial good stratification $\mathcal{V}$ of $X$ relative to $f.$

\begin{theorem}\label{Generalization 6.4 of DG}
Suppose that $g$ is tractable at the origin with respect to  $\mathcal{V}$ relative to $f.$ Then, for $0<|\delta|\ll\epsilon\ll1,$ 
\begin{center}
$B_{f,X}(0)-B_{f,X^g}(0)-\sum_{j=1}^{r}m_{f,b_j}(Eu_X(b_j)-Eu_{X^g}(b_j))=(-1)^{d-1}m,$
\end{center}
where $m$ is the number of stratified Morse critical points of a partial Morsification of \linebreak$g:X\cap f^{-1}(\delta)\cap B_{\epsilon}\rightarrow\mathbb{C}$ appearing on $X_{reg}\cap f^{-1}(\delta)\cap \{g\neq 0\}\cap B_{\epsilon}.$ 
\end{theorem} 
\noindent\textbf{Proof.} By Corollary \ref{Corollary 4.3 of DG}, if $0<|\delta|\ll\epsilon\ll 1,$\begin{center}
$\chi(X\cap f^{-1}(\delta)\cap B_{\epsilon},Eu_X)-\chi(X\cap g^{-1}(0)\cap f^{-1}(\delta)\cap B_{\epsilon},Eu_X)=(-1)^{d-1}m.$
\end{center}

If $V_i\nsubseteq \Sigma_{\mathcal{V}} g,$ $V_i$ intersects $\{g=0\}$ transversely and $Eu_X(V_i)= Eu_{X^g}(V_i\cap g^{-1}(0)).$ Let us denote by $W_1,\ldots, W_s$ the strata contained in $\Sigma_{\mathcal{V}}g.$ Then, 
\vspace{-0.3cm}
\begin{eqnarray*}
\chi(X\cap g^{-1}(0)\cap f^{-1}(\delta)\cap B_{\epsilon},Eu_X)&=&\sum_{V_i\nsubseteq \Sigma_{\mathcal{V}} g}\chi(V_i\cap g^{-1}(0)\cap f^{-1}(\delta)\cap B_{\epsilon})Eu_{X^g}(V_i\cap g^{-1}(0))\\
&+&\sum_{l=1}^s\chi(W_l\cap g^{-1}(0)\cap f^{-1}(\delta)\cap B_{\epsilon})Eu_X(W_l).
\end{eqnarray*}
For each $W_l\subseteq\Sigma_{\mathcal{V}} g,$ let $k_l$ be the number of branches $b_{l_t}$ containing in $W_l.$ Then, 
\begin{center}
$\chi(W_l\cap g^{-1}(0)\cap f^{-1}(\delta)\cap B_{\epsilon})=\displaystyle\sum_{b_{l_t}\subseteq W_l}\chi(b_{l_t}\cap f^{-1}(\delta)\cap B_{\epsilon})=\sum_{t=1}^{k_l} m_{f,b_{l_t}}$
\end{center}
and then
\begin{center}
$\displaystyle\sum_{l=1}^{s}\chi(W_l\cap g^{-1}(0)\cap f^{-1}(\delta)\cap B_{\epsilon})Eu_X(W_l)=\sum_{j=1}^{r}m_{f,b_j}Eu_X(b_j).$
\end{center}

On the other hand, 
\begin{eqnarray*}
B_{f,X^g}(0)&=&\sum_{V_i\nsubseteq \Sigma_{\mathcal{V}} g}\chi(V_i\cap g^{-1}(0)\cap f^{-1}(\delta)\cap B_{\epsilon})Eu_{X^g}(V_i\cap g^{-1}(0))\\
&+&\sum_{j=1}^rm_{f,b_j}Eu_{X^g}(b_j).
\end{eqnarray*}

Therefore, 
\begin{eqnarray}
\chi(X\cap g^{-1}(0)\cap f^{-1}(\delta)\cap B_{\epsilon},Eu_X)&=&B_{f,X^g}(0)\nonumber\\
&-&\sum_{j=1}^rm_{f,b_j}(Eu_{X^g}(b_j)-Eu_X(b_j)) \label{F31}
\end{eqnarray}

\fim




\section{Some results for functions with isolated singularity}

\hspace{0,5cm} In Section 6 of \cite{DG}, Dutertre and Grulha proved several relations between the Brasselet number of functions with isolated singularity and other invariants. In this section, we provide the generalization of some of their results to the context we describe in the following. Let $X$ be an analytic complex space and $\mathcal{W}=\{W_0, \ldots, W_q\}$ be a Whitney stratification of $X$ with $W_0=\{0\}.$ From now on, we consider $f$ and $g$ functions defined over $X$ such that $f$ has isolated singularity at the origin, $\Sigma_{\mathcal{W}}g$ is a one-dimensional analytic set and $\Sigma_{\mathcal{W}}g\cap\{f=0\}=\{0\}.$ Let $\mathcal{V}$ be the good stratification of $X$ induced by $f$ and suppose that $g$ is tractable at the origin with respect to $\mathcal{V}$ relative to $f.$ Notice that $\Sigma_{\mathcal{W}}g=\Sigma_{\mathcal{V}}g.$ 

\begin{lemma}\label{second stratification lemma}(\textbf{Second stratification lemma})
Let $\mathcal{V}$ be the good stratification of $X$ induced by $f, \mathcal{V}^{f}$ the collection of strata of $\mathcal{V}$ contained in $\{f=0\}$ and suppose that $g$ is tractable at the origin with respect to $\mathcal{V}$ relative to $f.$ Consider the refinement of $\mathcal{V},$\begin{eqnarray}\label{second stratification}
\mathcal{V}^{\prime\prime}=\{V_i\setminus\{g=0\}, V_i\cap\{g=0\}\setminus\Sigma_{\mathcal{W}}g, V_i\cap\Sigma_{\mathcal{W}}g, V_i\in\mathcal{V}\}\cup\{0\}.
\end{eqnarray}
Then $\mathcal{V}^{\prime\prime}$ is a good stratification of $X$ relative to $g$ such that $\mathcal{V}^{\prime\prime\{f=0\}},$ \begin{eqnarray*}
\mathcal{V}^{\prime\prime\{f=0\}}=\{V_i\cap \{f=0\}\setminus\{g=0\}, V_i\cap\{f=0\}\cap\{g=0\}\setminus\Sigma_{\mathcal{W}}g, V_i\in\mathcal{V}^{f}\}\cup\{0\},
\end{eqnarray*} 
\noindent is a good stratification of $X^f$ relative to $g|_{X^f}.$ 

Moreover, $f$ is prepolar at the origin with respect to $\mathcal{V}^{\prime\prime}$ relative to $g.$
\end{lemma}

\noindent\textbf{Proof.} Let us first show that $\mathcal{V}^{\prime\prime}$ is a good stratification of $X$ with respect to $g.$ 
\begin{enumerate}
    \item $V(g)$ is union of strata of type $V_i\cap\{g=0\}\setminus\Sigma_{\mathcal{W}}g$ and $V_i\cap\Sigma_{\mathcal{W}}g;$
    
    \item Let us show that $\{V_i\setminus\{g=0\},V_i\in\mathcal{V}\},$ which is equal to \begin{center}
        $\{W_i\cap\{f=0\}\setminus\{g=0\},W_i\setminus\big(\{f=0\}\cup\{g=0\}\big), W_i\in\mathcal{W}\},$ \end{center} is a Whitney stratification of $X\setminus\{g=0\}.$ Since $f$ has isolated singularity at the origin, \linebreak$\{f=0\}$ intersects each strata $W_i$ transversely. Therefore, since $\mathcal{W}$ is a Whitney stratification of $X,$ $\{V_i\setminus\{g=0\}, V_i\in\mathcal{V}\}$ satisfies Whitney's condition $(b)$.
    
    \item Let us verify the Thom condition. Let $V_{\lambda}\nsubseteq V(g), V_{\gamma}\subset V(g)$ be strata of $\mathcal{V}^{\prime\prime}$ and let $(p_k)$ be a sequence of points of $V_{\lambda}$ converging to a point $p\in V_{\gamma}.$ Suppose that the sequence of tangent spaces $T_{p_k}V(g|_{V_{\lambda}}-g|_{V_{\lambda}}(p_i))$ converges to $T.$ We must show that $T_p V_{\gamma}\subseteq T.$ If $p=0,$ then $V_{\gamma}=\{0\}$ and $\{0\}=T_p V_{\gamma}\subseteq T.$ 
    Suppose now that $p\neq 0$ and consider $V_{\gamma}=W_j\cap\Sigma_{\mathcal{W}}g, W_j\in\mathcal{W}.$ Since Thom stratifications always exist, one may take a refinement of $W_j\cap\Sigma_{\mathcal{W}}g$ that guarantees that the Thom condition is valid over this strata. Since $\Sigma_{\mathcal{W}}g$ is one-dimensional, this refinement would be given by taking off a finite number of points. Therefore, working on a sufficiently small neighborhood of the origin, Thom condition is verified over $W_j\cap\Sigma_{\mathcal{W}}g.$
    For $p\neq 0,$ we have two options for $V_{\lambda}\nsubseteq V(g),$ which are $W_i\setminus\big(\{f=0\}\cup\{g=0\}\big)$ and $W_i\cap\{f=0\}\setminus\{g=0\},$ where $W_i\in\mathcal{W}.$ 
    
    Suppose that $V_{\lambda}=W_i\setminus\big(\{f=0\}\cup\{g=0\}\big), W_i\in\mathcal{W},$ and let $\tilde{g}$ be an analytic extension of $g$ to an open neighborhood of the origin in $\mathbb{C}^n.$
    Then \begin{eqnarray*}
    \lim_{k\rightarrow\infty} T_{p_k}V(g|_{V_{\lambda}}-g|_{V_{\lambda}}(p_k))&=& \lim _{k\rightarrow\infty}T_{p_k}V(g|_{W_i\setminus\big(\{f=0\}\cup\{g=0\}\big)}-g|_{W_i\setminus\big(\{f=0\}\cup\{g=0\}\big)}(p_k))\\
    &=& \lim_{k\rightarrow\infty}T_{p_k} V(\tilde{g}-\tilde{g}(p_k))\cap T_{p_k} W_i\\
    &\subseteq&\lim_{k\rightarrow\infty}T_{p_k} V(\tilde{g}-\tilde{g}(p_k))\cap \lim_{k\rightarrow\infty} T_{p_k} W_i
    \end{eqnarray*} 
    
    If $p\in V_{\gamma}=V_j\cap\{g=0\}\setminus\Sigma_{\mathcal{W}}g, V_j\in\mathcal{V}$ and writing $\lim_{k\rightarrow\infty}T_{p_k} W_i=T_1,$ since $p\not\in\Sigma_{\mathcal{W}}g,$ the last limit is equal to \begin{eqnarray*}
    T_{p}V(\tilde{g}-\tilde{g}(p))\cap T_1= T_{p}V(\tilde{g})\cap T_1.
    \end{eqnarray*} 
    
    Suppose now that $V_{\gamma}=W_j\cap\{g=0\}\setminus\big(\{f=0\}\cup\Sigma_{\mathcal{W}}g\big), W_j\in\mathcal{W}.$ Then $$T_p V_{\gamma}=T_p\big(W_j\cap\{g=0\}\setminus(\{f=0\}\cup\Sigma_{\mathcal{W}}g)\big)=T_pW_j\cap T_pV(\tilde{g}).$$ By Whitney's condition $(a)$ over strata of $\mathcal{W}$, $T_{p}W_j\subseteq T_1.$ Since $g$ is tractable at the origin with respect to $\mathcal{V}$ relative to $f, T_p V(\tilde{g})$ intersects $T_p(W_j\setminus\{f=0\})=T_p W_j$ transversely at $p\not\in\Sigma_{\mathcal{W}}g.$ Therefore, the intersection $T_p V(\tilde{g})\cap T_1$ is transverse. Then we conclude that \begin{center}
        $\displaystyle\lim_{k\rightarrow\infty} T_{p_k}V(g|_{V_{\lambda}}-g|_{V_{\lambda}}(p_k))=T_p V(\tilde{g})\cap T_1.$
    \end{center} Therefore, $$T_pV_{\gamma}=T_pW_j\cap T_pV(\tilde{g})\subseteq T_1\cap T_pV(\tilde{g}).$$
    
    Now, let $\tilde{f}$ be an analytic extension of $f$ to the ambient space $U$ of $X.$ If \linebreak $V_{\gamma}=W_j\cap\{g=0\}\cap\{f=0\}\setminus\Sigma_{\mathcal{W}}g, W_j\in\mathcal{W},$ $$T_p V_{\gamma}=T_p(W_j\cap\{g=0\}\cap\{f=0\}\setminus\Sigma_{\mathcal{W}}g)=T_pW_j\cap T_pV(\tilde{g})\cap T_pV(\tilde{f}).$$ Using Whitney's condition $(a)$ over strata of $\mathcal{W}$ again, $$T_pV_{\gamma}=T_pW_j\cap T_pV(\tilde{g})\cap T_pV(\tilde{f})\subseteq T_1\cap T_pV(\tilde{g}).$$

    Let us now analyze the case where $V_{\lambda}=W_i\cap\{f=0\}\setminus\{g=0\},W_i\in\mathcal{W}.$ 
    
    Then \begin{eqnarray*}
    \lim_{k\rightarrow\infty} T_{p_k}V(g|_{V_{\lambda}}-g|_{V_{\lambda}}(p_k))&=& \lim _{k\rightarrow\infty}T_{p_k}V(g|_{W_i\cap\{f=0\}\setminus\{g=0\}}-g|_{W_i\cap\{f=0\}\setminus\{g=0\}}(p_k))\\
    &=& \lim_{k\rightarrow\infty}T_{p_k} V(\tilde{g}-\tilde{g}(p_k))\cap T_{p_k} W_i\cap T_{p_k} V(\tilde{f})\\
    &\subseteq&\lim_{k\rightarrow\infty}T_{p_k} V(\tilde{g}-\tilde{g}(p_k))\cap\lim_{k\rightarrow\infty} T_{p_k} W_i\cap \lim_{k\rightarrow\infty} T_{p_k} V(\tilde{f}).
    \end{eqnarray*} 
    
    Notice that $p$ must be contained in $\{f=0\},$ since it is a point of the closure of \linebreak $W_i\cap\{f=0\}\setminus\{g=0\}, W_i\in\mathcal{W}.$ Therefore, the only option we have for $V_{\gamma}$ is \linebreak $W_j\cap\{f=0\}\cap\{g=0\}\setminus\Sigma_{\mathcal{W}}g, W_j\in\mathcal{W}.$ Then, writing $\lim_{k\rightarrow\infty}T_{p_k} W_i=T_1,$ the last limit is equal to \begin{eqnarray*}
    T_{p}V(\tilde{g}-\tilde{g}(p))\cap T_1\cap T_p V(\tilde{f})= T_{p}V(\tilde{g})\cap T_1\cap T_p V(\tilde{f}).
    \end{eqnarray*} 
    
    For $V_{\gamma}=W_j\cap\{g=0\}\cap\{f=0\}\setminus\Sigma_{\mathcal{W}}g, W_j\in\mathcal{W},$ we have $$T_p V_{\gamma}=T_p(W_j\cap\{g=0\}\cap\{f=0\}\setminus\Sigma_{\mathcal{W}}g)=T_pW_j\cap T_pV(\tilde{g})\cap T_pV(\tilde{f}).$$
     Since $f$ has isolated singularity at the origin, $T_pV(\tilde{f})$ intersects $T_pW_j$ transversely and since $g$ is tractable at the origin with respect to $\mathcal{V}, T_pV(\tilde{g})$ intersects $T_pV(\tilde{f})\cap T_p W_j$ transversely.   Since, by Whitney's condition $(a)$ over strata of $\mathcal{W}, T_pW_j\subseteq T_1,$ the intersections on $T_{p}V(\tilde{g})\cap T_1\cap T_p V(\tilde{f})$ are transverse. Then we conclude that \begin{center}
         $ \lim_{k\rightarrow\infty} T_{p_k}V(g|_{V_{\lambda}}-g|_{V_{\lambda}}(p_k))=T_{p}V(\tilde{g})\cap T_1\cap T_p V(\tilde{f}).$
     \end{center} 
    Using Whitney's condition $(a)$ over strata of $\mathcal{W}$ again, $$T_pV_{\gamma}=T_pW_j\cap T_pV(\tilde{g})\cap T_pV(\tilde{f})\subseteq T_1\cap T_pV(\tilde{g})\cap T_p V(\tilde{f}).$$

\end{enumerate}

Let us now to verify that $\mathcal{V}^{\prime\prime\{f=0\}}$ is a good stratification of $X^f$ relative to $g|_{X^f}.$ This is valid because
\begin{center}
$ \mathcal{V}^{\prime\prime\{f=0\}}=\{W_i\cap\{f=0\}\setminus\{g=0\}, W_i\cap\{f=0\}\cap\{g=0\}\setminus\Sigma_{\mathcal{W}}g, W_i\in\mathcal{W}\}$
\end{center} is given by strata of $\mathcal{V}^{\prime\prime}.$ 

At last, we will show that $f$ is prepolar with respect to $\mathcal{V}^{\prime\prime}$ at the origin. For that, we need to verify that for all $V_{\alpha}\in\mathcal{V}^{\prime\prime}, 0\not\in V_{\alpha}, f|_{V_{\alpha}}$ is nonsingular. If $V_{\alpha}=V_i\setminus\{g=0\},V_i\in\mathcal{V},$ since $f$ has isolated singularity at the origin, $f|_{V_i\setminus\{g=0\}}$ has no singularity . Suppose now that $V_{\alpha}=V_i\cap\Sigma_{\mathcal{W}}g,$ with $V_i\in\mathcal{V}.$ Since, by Proposition 1.3 of \cite{Ms1}, $\Sigma f|_{V_{\alpha}}\subset\{f=0\}$ and, by hypothesis, $\Sigma_{\mathcal{W}}g\cap\{f=0\}=\{0\},$ $f|_{V_{\alpha}}$ is nonsingular.
Now, let $V_{\alpha}=V_i\cap\{g=0\}\setminus\Sigma_{\mathcal{W}}g, V_i\in\mathcal{V}$ and $x\in\Sigma f|_{V_{\alpha}}, x\neq0.$ Since $\Sigma f|_{V_{\alpha}}\subset\{f=0\}, V_i=W_i\cap\{f=0\},$ for some $W_i\in\mathcal{W}.$ Then $x\in W_i\cap\{f=0\}\cap\{g=0\}\setminus\Sigma_{\mathcal{W}}g.$ But $g$ is tractable at the origin with respect to $\mathcal{V}$ relative to $f$, which implies that $W_i\cap\{f=0\}$ intersects $\{g=0\}$ transversely and gives us a contradiction.  Therefore, $f$ is prepolar at the origin with respect to $\mathcal{V}^{\prime\prime}.$\fim



Let us see an adaptation of Theorem 3.9 of \cite{Ms1} to the case we are working on. 

\begin{lemma}\label{3.9 adapted}
Suppose that $g$ is tractable at the origin with respect to the good stratification $\mathcal{V}$ of $X$ induced by $f.$ Then, for $0<|\alpha|,|\delta|\ll\epsilon<1$ and a closed ball $B_{\epsilon}$ centered at the origin, \begin{center}
$\chi(X\cap g^{-1}(\alpha)\cap f^{-1}(\delta)\cap B_{\epsilon})=\chi(X\cap g^{-1}(\alpha)\cap f^{-1}(0)\cap B_{\epsilon}).$
\end{center}
\end{lemma} 
\noindent\textbf{Proof.} Let \begin{eqnarray}
\mathcal{V}^{\prime\prime}=\{V_i\setminus\{g=0\}, V_i\cap\{g=0\}\setminus\Sigma_{\mathcal{W}}g, V_i\cap\Sigma_{\mathcal{W}}g, V_i\in\mathcal{V}\}\cup\{0\}.
\end{eqnarray}

\noindent be the good stratification of $X$ relative to $g$ constructed in Lemma \ref{second stratification lemma}. By this lemma, $f$ is prepolar at the origin with respect to $\mathcal{V}^{\prime\prime}.$ So, $V(f)$ intersects each stratum of $\mathcal{V}^{\prime\prime}$ transversely in a neighborhood of the origin, except perhaps at the origin itself. Hence, we can choose a sufficiently small $\epsilon$ such that in an open ball containing $B_{\epsilon},$ $V(f)$ intersects \linebreak$\{V_{\lambda}\cap V(g)\setminus\Sigma_{\mathcal{W}} g, V_{\lambda}\in\mathcal{V}\}$ transversely and such that the sphere $\partial B_{\epsilon}$ intersects each $V_{\lambda}\cap V(g)\cap V(f)$ transversely.

Fixing the appropriate $\epsilon,$ let us show that, for $0<\eta,\nu\ll\epsilon,$ the map 

\begin{center}$\begin{array}{c}
B_{\epsilon}\cap X\cap\Phi^{-1}(int(D_{\eta})\times int(D_{\nu})-\Phi(\tilde{\Gamma}_{f,g}(\mathcal{V})\cup\Sigma_{\mathcal{W}} g))\\
\downarrow {}_{\Phi:=(f,g)}\\
int(D_{\eta})\times int(D_{\nu})-\Phi(\tilde{\Gamma}_{f,g}(\mathcal{V})\cup\Sigma_{\mathcal{W}} g)
\end{array}$
\end{center}
is a stratified proper submersion with respect to $\mathcal{V},$ where $D_{\eta}$ and $D_{\nu}$ are small closed balls centered at the origin. 

Since $X$ is Hausdorff, $B_{\epsilon}\cap X$ is compact and $(f,g)$ is a continuous map, $(f,g):B_{\epsilon}\cap X\rightarrow int(D_{\eta})\times int(D_{\nu})$ is proper and so it is the restriction $\Phi$ defined above.

Let us prove that $\Phi$ is a submersion. Since $f$ has isolated singularity at the origin and the symmetric relative polar curve $\tilde{\Gamma}_{f,g}(\mathcal{V})$ and the singular locus $\Sigma_{\mathcal{W}} g$ were excluded, $\Phi$ has no critical point inside $int(B_{\epsilon})\cap X.$ 

Let us now verify that $\Phi$ has no critical points on the boundary $\partial B_{\epsilon}\cap X.$ Let $\tilde{f}$ and $\tilde{g}$ be extensions of $f$ and $g$ to the ambient space, respectively. By contradiction, suppose that no matter how small we pick $\eta, \nu,$ $\Phi$ has a stratified critical point on the boundary $\partial B_{\epsilon}\cap X.$ Since the covering given by the stratification is locally finite, we can assume that all these critical points lie in some stratum $V_{\lambda}$. Then there exists a sequence  of critical points $(p_i)$ of $\partial B_{\epsilon}\cap V_{\lambda}$ such that $p_i\rightarrow p,$ $f(p_i)\rightarrow 0,\ g(p_i)\rightarrow 0,$ and \begin{equation}
T_{p_i}V(\tilde{f}-\tilde{f}(p_i))\cap T_{p_i}V(\tilde{g}-\tilde{g}(p_i))\cap T_{p_i}V_{\lambda}\subseteq T_{p_i}\partial B_{\epsilon}. \label{Contained in the sphere}
\end{equation}  

Hence, $p\in V_{\beta}\subset\overline{V_{\lambda}}$, $p\not\in\Sigma_{\mathcal{W}}g,$ $p\not\in\tilde{\Gamma}_{f,g}(\mathcal{V})$ and $f(p)=0.$ Then $V(\tilde{f})$, $V(\tilde{g})$ and $T_p V_{\beta}$ intersect transversely at $p$ and \begin{center}$T_{p_i}V(\tilde{f}-\tilde{f}(p_i))\rightarrow T_{p}V(\tilde{f}-\tilde{f}(0))= T_p V(\tilde{f})$ and $T_{p_i}V(\tilde{g}-\tilde{g}(p_i))\rightarrow T_{p}V(\tilde{g}-\tilde{g}(0))=T_p V(\tilde{g}).$\end{center}

Therefore, if we suppose that $T_{p_i}V_{\lambda}\rightarrow \mathcal{T},$ applying the limit to (\ref{Contained in the sphere}) and using Whitney's condition $(a),$ we obtain that \begin{center}
$T_p V(\tilde{f})\cap T_p V(\tilde{g})\cap T_{p}V_{\beta}\subseteq T_p V(\tilde{f})\cap T_p V(\tilde{g})\cap \mathcal{T}\subseteq T_{p}\partial B_{\epsilon},$
\end{center}
\noindent which is a contradiction, since we choose $\epsilon$ sufficiently small such that \linebreak$V(g)\cap V(f)\cap V_{\beta}$ intersects $\partial B_{\epsilon}$ transversely. 

Hence, $\Phi$ is a stratified proper submersion. By Ehresmann Fibration Theorem, all fibres are homeomorphic. 

Notice that, since $g$ is tractable at the origin with respect to $\mathcal{V}$ relative to $f,$ $\Phi$ has no critical points contained in $V(f),$ that is, $\Phi$ has no critical points of the type $(0,\alpha),\alpha\neq0.$ 

Then, for $0<|\alpha|\ll|\delta|\ll\epsilon<1,$ with $\alpha$ being a regular value of $g,$ the fibres $\Phi^{-1}(\delta,\alpha)$ and $\Phi^{-1}(0,\alpha)$ are homeomorphic, that is, \begin{center}
$\chi(X\cap g^{-1}(\alpha)\cap f^{-1}(\delta)\cap B_{\epsilon})=\chi(X\cap g^{-1}(\alpha)\cap f^{-1}(0)\cap B_{\epsilon}).$
\end{center}
\fim

Another property we will need is the following.

\begin{lemma}\label{critical points intersection}
Let $\mathcal{V}$ be the good stratification of $X$ induced by $f$ and suppose that $g$ is tractable at the origin with respect to $\mathcal{V}$ relative to $f$. If $0<|\delta|\ll\epsilon\ll 1,$ then \begin{center}$\Sigma_{\mathcal{W}\cap\{f=\delta\}} g \cap\{g=0\}\cap B_{\epsilon}=\Sigma_{\mathcal{W}} g\cap\{f=\delta\}\cap \{g=0\}\cap B_{\epsilon}.$\end{center}
\end{lemma}

\noindent\textbf{Proof.} Let $\tilde{g}$ and $\tilde{f}$ be analytic extensions of $g$ and $f$ to the ambient space $U.$ 

Let $p\in\Sigma_{\mathcal{W}} g\cap\{f=\delta\}\cap \{g=0\}\cap B_{\epsilon}$ and $V_{\alpha}$ be the stratum of $\mathcal{V}$ that contains $p.$ Then $d_p\tilde{g}|_{V_{\alpha}}=0,$ and $rk(d_p\tilde{g}|_{V_{\alpha}},d_p\tilde{f}|_{V_{\alpha}})\leq1.$ So, $p$ is a critical point of $g|_{V_{\alpha}\cap\{f=\delta\}},$ that is, \linebreak$p\in \Sigma_{\mathcal{W}\cap \{f=\delta\}} g|\cap \{g=0\}\cap B_{\epsilon}.$

Suppose now that there exists in $\Sigma_{\mathcal{W}\cap \{f=f(p_i)\}} g\cap\{g=0\}\cap B_{\epsilon}\setminus\Sigma_{\mathcal{W}}g$ a sequence of points $(p_i)$ converging to $0.$ Then, for all $i, p_i\in\Gamma_{f,g}(\mathcal{V})\cap\{g=0\}\setminus\Sigma_{\mathcal{W}}g.$ Since $(p_i)\in\{g=0\}\setminus\Sigma_{\mathcal{W}}g,$ each $p_i$ is a critical point of $f|_{\{g=0\}\setminus\Sigma_{\mathcal{W}}g}.$ So, by Proposition 1.3 of \cite{Ms1}, $p_i\in\{f=0\},$ for all $i,$ which is a contradiction.\fim

If $V_1,\ldots,V_q$ are the strata not contained in $\{f=0\},$ we can write $\Sigma_{\mathcal{W}}g=\linebreak b_1\cup\ldots\cup b_r$ as a union of branches $b_j,$ where $b_j\subseteq V_{i(j)},$ for some $i(j)\in\{1,\ldots,q\},$ as we saw before. Let $\delta$ be a regular value of $f, 0<|\delta|\ll1,$ and  $f^{-1}(\delta)\cap b_j=\{x_{i_1},\ldots,x_{i_{k(j)}}\}.$ For each $x_{l}\in f^{-1}(\delta)\cap b_j,$ let $D_{x_l}$ be the closed ball with center at $x_{l}$ and radius $0<r_{l}\ll1.$ We choose $r_{l}$ sufficiently small such that the balls $D_{x_l}$ are pairwise disjoint and the union of balls $D_j=D_{x_{i_1}}\cup\ldots\cup D_{x_{i_{k(j)}}}$ is contained in $B_{\epsilon},$ where $0<|\delta|\ll\epsilon\ll1$ and $\epsilon$ is sufficiently small such that the local Euler obstruction of $X$ is constant on $b_j\cap B_{\epsilon}$. Notice that, in this case, we can choose $x_{l}\in b_j, j\in\{1,\ldots,r\}, l\in\{i_1,\ldots, i_{k(j)}\},$ and write $Eu_X(x_{l})=Eu_X(b_j).$

Before we prove the first theorem of this section, we will prove a useful regularity condition over the branches $b_j.$  Notice that the next lemma is a version of Corollary 2.2 of \cite{iomdin1974complex} in our setting.

\begin{lemma}\label{fibres lemma}
Let $0<|\alpha|\ll1$ and $0<|\delta|\ll1$ be regular values of $g$ and $f,$ respectively.  For all $l_1\neq l_2,l_1,l_2\in\{i_1,\ldots, i_{k(j)}\},$ and $|\alpha|\ll|\delta|,$ \begin{center}
$\chi(X\cap g^{-1}(\alpha)\cap f^{-1}(\delta)\cap D_{x_{l_1}})=\chi(X\cap g^{-1}(\alpha)\cap f^{-1}(\delta)\cap D_{x_{l_2}}).$
\end{center}
\end{lemma}
\noindent\textbf{Proof.} Consider the function $\varphi$ over $b_j\setminus\{0\}$ given by $\varphi(x)=\chi(X\cap\{g=\alpha\}\cap\{f=f(x)\}\cap D_{x}),$ where $|\alpha|\ll|f(x)|$ and $D_{x}$ the closed ball with center at $x$ and radius $0<r_{x}\ll1.$ We should prove that $\varphi$ is constant. Since $b_j\setminus\{0\}$ is connected, it is sufficient to show that $\varphi$ is locally constant, that is, given $x\in b_j\setminus\{0\},$ it must exist a neighborhood $V_x$ of $x$ such that for all $y\in b_j\setminus\{0\}\cap V_x, \ \varphi(x)=\varphi(y).$ It is enough to show that there exist $\epsilon_x>0$ and a neighborhood $V_x$ such that for all $y\in b_j\setminus\{0\}\cap V_x$ and all $0<\epsilon\leq\epsilon_x, S(y,\epsilon)$ intersects $g^{-1}(0)\cap f^{-1}(f(y))$ transversely.

Let us denote by $N(\epsilon)$ the tube $\{z\in U; d(\Sigma_{\mathcal{V}}g,z)=\epsilon\}.$ We can replace $S(y,\epsilon)$ with $N(\epsilon)$ and we have to show that there exist $\epsilon_x>0$ and a neighborhood $V_x$ of $x$ such that for all $y\in b_j\setminus\{0\}\cap V_x$ and all $\epsilon\leq\epsilon_x,$ $N(\epsilon)$ intersects $g^{-1}(0)\cap f^{-1}(f(y))$ transversely. We can also replace the distance function to  $\Sigma_{\mathcal{W}} g$ with a real analytic function $h$ such that $h^{-1}(0)=\Sigma_{\mathcal{W}} g$ and $h\geq0.$ Then we replace the tube $N(\epsilon)$ with $\{z\in U;h(z)=\epsilon\}=h^{-1}(\epsilon).$

By contradiction, suppose  that there exists a point $y$ such that $N(\epsilon)$ does not intersect $g^{-1}(0)\cap f^{-1}(f(y))$ transversely. Let $V_{\alpha}$ be the stratum of $\mathcal{V}$ that contains $y.$ Then, using the terminology of Iomdin in \cite{iomdin1974complex}, since $y\in V_{\alpha}\cap\{g=0\},$ the vectors $grad \ h|_{V_{\alpha}}(y)$ and $grad \ f|_{V_{\alpha}}(y)=(1/2 f|_{V_{\alpha}}(y))grad \ ||f|_{V_{\alpha}}(y)||$ are complex linearly dependent $mod \ grad \ g|_{V_{\alpha}}(y).$ Hence, $grad \ f|_{V_{\alpha}}(y)=\lambda grad \ h|_{V_{\alpha}}(y), mod \ grad \ g|_{V_{\alpha}}(y).$



Now, \begin{eqnarray*}
    grad \ ||f|_{V_{\alpha}}||^2(y)&=&2f|_{V_{\alpha}}(y)grad \ f|_{V_{\alpha}}(y)=2f|_{V_{\alpha}}(y)\lambda grad \ h|_{V_{\alpha}}(y)\\
    &=&\frac{\lambda f|_{V_{\alpha}}(y)}{h|_{V_{\alpha}}(y)}2h|_{V_{\alpha}}(y)grad \ h|_{V_{\alpha}}(y)=\gamma grad \ ||h|_{V_{\alpha}}||^2(y), 
\end{eqnarray*}
\noindent with $\gamma=\frac{\lambda f|_{V_{\alpha}}(y)}{h|_{V_{\alpha}}(y)}.$
The last equality means that the vectors $grad \ ||f|_{V_{\alpha}}||^2(y)$ and $grad \ ||h|_{V_{\alpha}}||^2(y)$ are complex linearly dependent $mod \  grad \ g|_{V_{\alpha}}(y).$ This contradicts Corollary 1.7 of \cite{iomdin1974local} using the functions $||h|_{V_{\alpha}}||^2$ and $||f|_{V_{\alpha}}||^2,$ since $\Sigma_{\mathcal{W}} g\cap\{f=0\}=\{0\},$ this corollary implies that there exist $\epsilon>0$ and a neighborhood $G$ of $\Sigma_{\mathcal{W}} g$ in $\{g=0\}$ such that at points $z$ of $D_{\epsilon}\cap G\setminus\Sigma_{\mathcal{W}} g, grad \ ||f|_{V_{\alpha}}||^2(z)$ and $grad \ ||h|_{V_{\alpha}}||^2(z)$ are complex linearly independent $mod \ grad \ g|_{V_{\alpha}}(z).$   

Therefore, $\varphi$ is locally constant.
\fim

\begin{observation}\label{independence of delta}
The last lemma shows that, for $0\leq|\alpha|\ll|\delta|\ll\epsilon\ll1,$ the Euler characteristic of $X\cap g^{-1}(\alpha)\cap f^{-1}(\delta)\cap D_{x_l}$ is constant over $b_j\cap B_{\epsilon}, j\in\{1,\ldots,r\}$ and $l\in\{i_1\ldots,i_{k(j)}\}.$ Then, for each stratum $V_1,\ldots,V_q$ of $\mathcal{V}$ not contained in $\{f=0\}, \chi(\overline{V_i}\cap g^{-1}(\alpha)\cap f^{-1}(\delta)\cap D_{x_l})$ is constant over $b_j\cap B_{\epsilon}.$ Notice that $\chi(V_i\cap g^{-1}(\alpha)\cap f^{-1}(\delta)\cap D_{x_l})$ is also constant over $b_j\cap B_{\epsilon}.$ If $V_i\neq\overline{V_i},$ write $\overline{V_i}=V_i\cup (\overline{V_i}\setminus V_i),$ use that $\overline{V_i}\setminus V_i$ is analytic and union of strata of $\mathcal{V}$ of smaller dimension and notice that for $V_0=\{0\},$ we have $\overline{V_0}=V_0.$ 
\end{observation} 

Let $\beta:X\rightarrow\mathbb{Z}$ be a constructible function with respect to $\mathcal{W}.$ Using Remark \ref{independence of delta}, since each $b_j$ is contained in one unique stratum of $\mathcal{W}$ and $\beta$ is constant over each one of them, we can use the following notation: \begin{enumerate}
    \item $\beta(b_j):=\beta(x_{l}),$ for a chosen $x_{l}\in b_j;$
    \item $\tilde{\beta}(b_j)=\chi(X\cap g^{-1}(\alpha)\cap f^{-1}(\delta)\cap D_{x_{l}},\beta)$ for $j\in\{1,\ldots,r\}$ and $l\in\{i_1,\ldots,i_{k(j)}\}.$
\end{enumerate}

\begin{theorem}\label{Euler characteristic 1}
Let $\beta:X\rightarrow\mathbb{Z}$ be a constructible function with respect to the stratification $\mathcal{W}$. For $0<|\alpha|\ll |\delta|\ll\epsilon \ll 1$, we have\begin{eqnarray*}
\chi(X\cap g^{-1}(\alpha)\cap f^{-1}(0)\cap B_{\epsilon},\beta)&=&\chi(X\cap g^{-1}(0)\cap f^{-1}(\delta)\cap B_{\epsilon},\beta)\\
&-&\sum_{j=1}^{r}m_{f,b_j}(\beta(b_j)-\tilde{\beta}(b_j)).
\end{eqnarray*}
\end{theorem}

\noindent\textbf{Proof.} We have
\begin{eqnarray*}
\chi(X\cap g^{-1}(0)\cap f^{-1}(\delta)\cap B_{\epsilon})
&=& \chi(X\cap g^{-1}(0)\cap f^{-1}(\delta)\cap B_{\epsilon}\setminus \cup_{j=1}^{r}D_j)\\
&+&\sum_{j=1}^{r}\chi(X\cap g^{-1}(0)\cap f^{-1}(\delta)\cap D_j).
\end{eqnarray*}

For each $j\in\{1,\ldots,r\},$ as we saw before, $D_j=D_{x_{i_1}}\cup\ldots D_{x_{i_{k(j)}}},$ where $D_{x_l}$ is a closed ball with center at $x_{l}.$ Since $f$ is an analytic function germ and, for each $l\in\{i_1,\ldots,i_{k(j)}\},$\linebreak $X\cap g^{-1}(0)\cap f^{-1}(\delta)\cap D_{x_l}$ is an analytic germ at $x_{l},$ it is contractible. So,\linebreak $\chi(X\cap g^{-1}(0)\cap f^{-1}(\delta)\cap D_{x_l})=1.$ Hence, $\sum_{j=1}^{r}\chi(X\cap g^{-1}(0)\cap f^{-1}(\delta)\cap D_j)=\sum^{r}_{j=1}m_{f,b_j}.$

On the other hand, by Lemma \ref{fibres lemma}, fixing $l\in\{i_1,\ldots,i_{k(j)}\},$
\begin{eqnarray*}
\chi(X\cap g^{-1}(\alpha)\cap f^{-1}(\delta)\cap B_{\epsilon})&=&\chi(X\cap g^{-1}(\alpha)\cap f^{-1}(\delta)\cap B_{\epsilon}\setminus \cup_{j=1}^{r}D_j)\\&+&\sum_{j=1}^{r}m_{f,b_j}\chi(X\cap g^{-1}(\alpha)\cap f^{-1}(\delta)\cap D_{x_l}).
\end{eqnarray*}

By Lemma \ref{critical points intersection}, $\Sigma_{\mathcal{W}} g\cap\{f=\delta\}\cap\{g=0\}\cap B_{\epsilon}=\{x_1,\ldots,x_s\}$ is the set of critical points of $g|_{\{f=\delta\}\cap B_{\epsilon}}$ appearing in $\{g=0\}.$ By Lemma \ref{3.9 adapted} and for $0<|\alpha|\ll |\delta|\ll\epsilon \ll 1,$ 
 $$\chi(X\cap g^{-1}(\alpha)\cap f^{-1}(0)\cap B_{\epsilon})=\chi(X\cap g^{-1}(\alpha)\cap f^{-1}(\delta)\cap B_{\epsilon})$$
\vspace{-0,5cm}
\begin{eqnarray*}
&=&\chi(X\cap g^{-1}(\alpha)\cap f^{-1}(\delta)\cap B_{\epsilon}\setminus (\cup_{j=1}^{r} D_{j}))+\sum_{j=1}^{r}m_{f,b_j}\chi(X\cap g^{-1}(\alpha)\cap f^{-1}(\delta)\cap D_{x_l})\\ &=& \chi(X\cap g^{-1}(0)\cap f^{-1}(\delta)\cap B_{\epsilon}\setminus (\cup_{j=1}^{r} D_j))+\sum_{j=1}^{r}m_{f,b_j}\chi(X\cap g^{-1}(\alpha)\cap f^{-1}(\delta)\cap D_{x_l})\\ &=& \chi(X\cap g^{-1}(0)\cap f^{-1}(\delta)\cap B_{\epsilon})-\sum_{j=1}^{r} m_{f,b_j}(1- \chi(X\cap g^{-1}(\alpha)\cap f^{-1}(\delta)\cap D_{x_l})).
\end{eqnarray*}

By additivity of the constructible function $\beta$, we obtain the formula.\fim


\begin{observation}
Let $\mathcal{W}$ be a Whitney stratification of $X$ and $\mathcal{V}$ be the good stratification of $X$ induced by $f.$ Suppose that $g$ is tractable at the origin with respect to $\mathcal{V}$ relative to $f,$ $ \Sigma_{\mathcal{W}}g$ is one-dimensional and that $\Sigma_{\mathcal{W}}g\cap\{f=0\}=\{0\}.$ The refinement $\mathcal{V}^{\prime\prime}$ of $\mathcal{V}^{\prime},$ constructed in Lemma \ref{second stratification lemma}, is a refinement of the stratification constructed in Lemma \ref{first stratification lemma}. Therefore, we can refine a Whitney stratification of $X$ to obtain an appropriate stratification for which the Brasselet numbers $B_{f,X}(0), B_{f,X^g}(0), B_{g,X(0)}$ and $B_{g,X^f}(0)$ can be explicitly calculated.
\end{observation}

Applying the previous theorem to the case where $\beta=Eu_X,$ we can compare $B_{g,X^f}(0)$ and $B_{f,X^g}(0).$ 

\begin{corollary}\label{Brasselet number 1}
Suppose that $g$ is tractable at the origin with respect to $\mathcal{V}$ relative to $f.$ Then, for $0\ll|\delta|\ll\epsilon\ll1,$ \begin{center}
$B_{g,X^f}(0)=B_{f,X^g}(0)-\sum_{j=1}^{r}m_{f,b_j}(Eu_{X^g}(b_j)-B_{g,X\cap f^{-1}(\delta)}(b_j)).$
\end{center}
\end{corollary}

\noindent\textbf{Proof.} Applying Theorem \ref{Euler characteristic 1} to $\beta=Eu_X$, we obtain \begin{eqnarray*}
\chi(X\cap g^{-1}(\alpha)\cap f^{-1}(0)\cap B_{\epsilon},Eu_X)&=&\chi(X\cap g^{-1}(0)\cap f^{-1}(\delta)\cap B_{\epsilon},Eu_X)\\
&-&\sum_{j=1}^{r}m_{f,b_j}(Eu_X(b_j)-\tilde{Eu}_X(b_j)).
\end{eqnarray*}

Let $W_1,\ldots,W_t$ be the strata of $\mathcal{V}^{\prime\prime}$ not contained in $\{g=0\}.$ Since $g$ is tractable at the origin with respect to $\mathcal{V}^{\prime}$ relative to $f,$ by Lemma \ref{second stratification lemma}, $f$ is prepolar a the origin with respect to $\mathcal{V}^{\prime\prime},$ that is, $\{f=0\}$ intersects each $W_i$ transversely, for $i\in\{1,\ldots,t\}.$ So, $Eu_X(W_i)=Eu_{X^f}(S),$ for each connected component $S$ of $W_i^f.$ Then, for $0<|\alpha|\ll\epsilon \ll 1,$ \begin{eqnarray*}
\chi(X\cap f^{-1}(0)\cap g^{-1}(\alpha)\cap B_{\epsilon},Eu_X)&=&\sum_{i=1}^{t}\chi(W_i\cap f^{-1}(0)\cap g^{-1}(\alpha)\cap B_{\epsilon})Eu_X(W_i)\\&=&\sum_{i=1}^{t}\sum_S\chi(S\cap f^{-1}(0)\cap g^{-1}(\alpha)\cap B_{\epsilon})Eu_{X^f}(S)\\&=&B_{g,X^f}(0).
\end{eqnarray*}

Using Equation (\ref{F31}) of Theorem \ref{Generalization 6.4 of DG}, we have that \begin{center}
$\chi(X\cap g^{-1}(0)\cap f^{-1}(\delta)\cap B_{\epsilon},Eu_X)=B_{f,X^g}(0)+\sum_{j=1}^{r}m_{f,b_j}(Eu_X(b_j)-Eu_{X^g}(b_j)).$
\end{center}

Consider the strata of $\mathcal{V}^{\prime\prime}$ not contained in $\{f=0\}$ or in $\{g=0\},$  \begin{center}
    $\{V_1\setminus\{f=0\}\cup\{g=0\},\ldots,V_q\setminus\{f=0\}\cup\{g=0\}, V_i\in\mathcal{V}\}.$
    
\end{center} Since $f$ is prepolar at the origin with respect to $\mathcal{V}^{\prime\prime}$, $f^{-1}(\delta)$ intersects each stratum \linebreak $U_j=V_j\setminus\{f=0\}\cup\{g=0\}$ transversely. Then \begin{eqnarray*}
\chi(X\cap g^{-1}(\alpha)\cap f^{-1}(\delta)\cap D_{x_l},Eu_X)&=&\sum_{i=1}^{q}\chi(U_i\cap g^{-1}(\alpha)\cap f^{-1}(\delta)\cap D_{x_l})Eu_X(U_i)\\&=&
\sum^q_{i=1}\chi(U_i\cap g^{-1}(\alpha)\cap f^{-1}(\delta)\cap D_{x_l})Eu_{X\cap f^{-1}(\delta)}(U_i\cap f^{-1}(\delta))\\&=&B_{g,X\cap f^{-1}(\delta)}(x_{l})\\&=&B_{g,X\cap f^{-1}(\delta)}(b_j),
\end{eqnarray*}

\noindent where the last equality holds by Remark \ref{independence of delta}.\fim


  Let $\mathcal{W}=\{\{0\}, W_1\ldots, W_q\}$ be a Whitney stratification of $X$, $\mathcal{V}$ be the good stratification of $X$ induced by $f, \mathcal{V}^{\prime}$ the good stratification of $X$ relative to $f$ obtained as a refinement of $\mathcal{V}$ in Lemma \ref{first stratification lemma} and $\mathcal{V}^{\prime\prime}$ be the good stratification of $X$ relative to $g$ obtained as a refinement of $\mathcal{V}$ in Lemma \ref{second stratification lemma}. Suppose that $\Sigma_{\mathcal{W}}g\cap \{f=0\}=\{0\}.$ Let $T_1,\ldots,T_{q}$ be the strata of $\mathcal{V}^{\prime\prime}$ not contained in $\{g=0\}$ and  $V_1,\ldots,V_{q}$ be the strata of $\mathcal{V}^{\prime}$ not contained in $\{f=0\}.$ Let $n_s$ (resp. $m_t$) be the number of stratified Morse critical points of a Morsification of \linebreak$f:X\cap g^{-1}(\alpha)\cap B_{\epsilon}\rightarrow\mathbb{C}$ (resp. $g:X\cap f^{-1}(\delta)\cap B_{\epsilon}\rightarrow\mathbb{C}$) appearing on \linebreak $T_s\cap g^{-1}(\alpha)\cap \{f\neq 0\}\cap B_{\epsilon}$ (resp. $V_t\cap f^{-1}(\delta)\cap \{g\neq 0\}\cap B_{\epsilon}$), where $0<|\delta|\ll1$ is a regular value of $f$ and $0<|\alpha|\ll1$ is a regular value of $g.$

\begin{theorem}\label{Euler Characteristic 2}
Let $\beta:X\rightarrow\mathbb{Z}$ be a constructible function with respect to $\mathcal{W}$ and suppose that $g$ is tractable at the origin with respect to $\mathcal{V}$ relative to $f$. For $0<\mid\alpha\mid\ll\mid\delta\mid\ll\epsilon\ll 1,$ \begin{eqnarray}
\chi(X\cap g^{-1}(\alpha)\cap B_{\epsilon},\beta)&-&\chi(X\cap f^{-1}(\delta)\cap B_{\epsilon},\beta)=\sum_{s=0}^{q}(-1)^{\dim T_s-1}n_s\eta(W_s,\beta)\nonumber\\
&-&\sum_{t=0}^{q}(-1)^{\dim V_t-1}m_{t}\eta(V_{t},\beta)-\sum_{j=1}^{r}m_{f,b_j}(\beta(b_j)-\tilde{\beta}(b_j)).\nonumber
\end{eqnarray}
\end{theorem}

\noindent\textbf{Proof.} By Lemma \ref{second stratification lemma}, since $g$ is tractable at the origin with respect to $\mathcal{V}^{\prime}$ relative to $f$ and $\Sigma_{\mathcal{W}}g\cap \{f=0\}=\{0\},$ $f$ is prepolar at the origin with respect to $\mathcal{V}^{\prime\prime}$ and, therefore, tractable at the origin with respect to $\mathcal{V}^{\prime\prime}$ relative to $g.$   
By Theorem \ref{4.2 DG}, \begin{eqnarray*}
\chi(X\cap g^{-1}(\alpha)\cap B_{\epsilon},\beta)-\chi(X\cap g^{-1}(\alpha)\cap f^{-1}(0)\cap B_{\epsilon},\beta)=\sum_{s=0}^{q}(-1)^{\dim T_s-1}n_s\eta(T_s,\beta).
\end{eqnarray*}

Since $g$ is tractable at the origin with respect to $\mathcal{V}$ relative to $f,$ also by Theorem \ref{4.2 DG}, \begin{eqnarray*}
\chi(X\cap f^{-1}(\delta)\cap B_{\epsilon},\beta)-\chi(X\cap g^{-1}(0)\cap f^{-1}(\delta)\cap B_{\epsilon},\beta)&=&\sum_{t=0}^{q}(-1)^{\dim V_t-1}m_{t}\eta(V_{t},\beta).
\end{eqnarray*} 

Using Theorem \ref{Euler characteristic 1}, we have the formula. \fim


If we apply Theorem \ref{Euler Characteristic 2} to the case where $\beta=Eu_X$, we obtain the following consequence.

\begin{corollary}\label{generalization 6.5}
Suppose that $g$ is tractable at the origin with respect to the good stratification $\mathcal{V}$ of $X$ induced by $f.$ For $0<|\alpha|\ll|\delta|\ll\epsilon\ll1,$
 \begin{center}
$B_{g,X}(0)-B_{f,X}(0)=(-1)^{d-1}(n_{reg}-m_{reg})-\sum_{j=1}^{r}m_{f,b_j}(Eu_X(b_j)-B_{g,X\cap\{f=\delta\}}(b_j)),$
\end{center}
\noindent where $n_{reg}=n_q$ and $m_{reg}=m_{q}$ in the previous notation.
\end{corollary}

\noindent\textbf{Proof.} First we have $\eta(T_s,Eu_X)=0,$ for $s\in \{1,\ldots,q-1\},$  $\eta(V_t,Eu_X)=0,$ for \linebreak$t\in \{1,\ldots,q-1\},$ where $V_t\in\mathcal{V}^{\prime\prime}$ is the strata not contained in $\{f=0\}$ and $W_s\in\mathcal{V}^{\prime\prime}$ is the strata not contained in $\{g=0\}.$ Also, since the local Euler obstruction is constant over $b_j\cap B_{\epsilon},$ we can write $Eu_X(x_{l})=Eu_X(b_j)$ and, by Lemma \ref{fibres lemma},  $B_{g,X\cap\{f=\delta\}}(x_{l})=B_{g,X\cap\{f=\delta\}}(b_j), l\in\{i_1,\ldots, i_{k(j)}\}.$ \fim



\vspace{0.1 cm}







\section{Applications to generic linear forms}

\hspace{0,5cm} In this section, we apply the results of the previous sections taking for $f$ the restriction to $X$ of a generic linear form $l:\mathbb{C}^n\rightarrow\mathbb{C}.$

\begin{lemma}\label{generic linear form lemma}
Let $V\subset\mathbb{C}^N$ be an analytic complex subset of dimension $d$ and $l:\mathbb{C}^N\rightarrow\mathbb{C}$ be a generic linear form. Then $l^{-1}(0)$ is transverse to $V\setminus\{0\}.$
\end{lemma}

\noindent\textbf{Proof.} We fix local coordinates $(x_1,\ldots,x_N)$ in $\mathbb{C}^N$ and define, for each $a=(a_1,\ldots,a_N)\in\mathbb{C}^N,$ $l_a(x)=a_1x_1+\cdots+a_Nx_N.$ Let $W=\{(x,a)\in\mathbb{C}^N\times\mathbb{C}^N;x\in V\setminus\{0\}, l_a(x)=0\}.$ Then $\dim W=2N-(N-d+1)=N+d-1.$ Consider the projection $\pi:W\rightarrow\mathbb{C}^N$ given by \linebreak $(x,a)\mapsto a$ and let $\Delta\subset\mathbb{C}^N$ be the discriminant of $\pi.$ 

If $a\in\mathbb{C}^N\setminus\Delta,$ then $V\setminus\{0\}$ intersects $\{l_a=0\}$ transversely and $\pi^{-1}(a),$ which is equal to $V\setminus\{0\}\cap\{l_a|_{V\setminus\{0\}}=0\},$ is a submanifold of $W$ with dimension $\dim W-N=N+d-1-N=d-1.$ \fim

Let $\mathcal{W}=\{\{0\}, W_1\ldots,W_q\}$ be a Whitney stratification of $X$ and $\mathcal{V}$ be the good stratification induced by a generic linear form $l.$ Then Lemma \ref{first stratification lemma} can be applied to $\mathcal{V}$ and we obtain a good stratification $\mathcal{V}^{\prime}$ of $X$ relative to $l$ such that $\mathcal{V}^{\prime\{g=0\}}$ is a good stratification of $X^g$ relative to $l_{|\{g=0\}}.$ Also, Lemma \ref{second stratification lemma} provides a good stratification $\mathcal{V}^{\prime\prime}$ of $X$ relative to $g$ constructed as a refinement o $\mathcal{V}.$

\begin{lemma}\label{generacity of being tractable}
For a generic linear form $l,$ the function $g$ is tractable at the origin with respect to the good stratification $\mathcal{V}$ of $X$ induced by $l$.
\end{lemma}
\noindent\textbf{Proof.} Let $\mathcal{V}=\{W_i\setminus\{l=0\}, W_i\cap\{l=0\}, W_i\in\mathcal{W}\}$ be the good stratification of $X$ induced by $l.$ By \cite{Ms1}, $l$ is tractable at the origin with respect to the good stratification $\mathcal{V}^{\prime\prime}$ relative to $g,$ that is, $\dim\tilde{\Gamma}_{l,g}(\mathcal{V}^{\prime\prime})\leq 1$ and $l|_{V_{\alpha}},$ with $V_{\alpha}\in X^g,$ has no critical points in a neighborhood of the origin except perhaps the origin itself. Now, \begin{eqnarray*}
\tilde{\Gamma}_{l,g}(\mathcal{V}^{\prime\prime})=\cup_{V_i\in\mathcal{V}^{\prime\prime}}\overline{\Sigma(l,g)|_{V_i\setminus(\{l=0\}\cup\{g=0\})}}=\cup_{W_i\in\mathcal{W}}\overline{\Sigma(l,g)|_{W_i\setminus(\{l=0\}\cup\{g=0\})}}=\tilde{\Gamma}_{l,g}(\mathcal{V}).
\end{eqnarray*}
Hence, $\dim\tilde{\Gamma}_{l,g}(\mathcal{V})\leq 1.$ Consider the strata $W_i\cap\{l=0\}$ in $\mathcal{V}.$ Let $G$ and $L$ be analytic extensions of $g$ and of $l$ in a neighborhood of the origin, respectively, and $x$ be a critical point of $g$ in $W_i\cap\{l=0\}.$ Then $d_xG|_{W_i\cap\{l=0\}}=\lambda(x)d_xL|_{W_i\cap\{l=0\}}.$ Since $x\in\{l=0\}$ and $x\neq0, x\not\in\Sigma_{\mathcal{W}}g.$ Therefore, $\lambda(x)\neq0$ and $d_xL|_{W_i\cap\{l=0\}}=\frac{1}{\lambda(x)}d_xG|_{W_i\cap\{l=0\}}.$ Hence, $x$ is a critical point of $l|_{W_i\cap\{g=0\}\cap\{l=0\}\setminus\Sigma_{\mathcal{W}}g},$ which is a contradiction.\fim

 A consequence of Theorem \ref{Generalization 6.4 of DG} is a relation between the differences of the Euler obstruction at the origin and at the branches.  
\begin{corollary}\label{Difference Euler obstructions}
Let $m$ be the number of stratified Morse critical points of a partial Morsification of $g:X\cap l^{-1}(\delta)\cap B_{\epsilon}\rightarrow\mathbb{C}$ appearing on $X_{reg}\cap l^{-1}(\delta)\cap \{g\neq 0\}\cap B_{\epsilon}$ and $m_{b_j}$ be the multiplicity of the branch $b_j$ at the origin. For $0<|\delta|\ll\epsilon\ll1,$ we have,
 
\begin{center}
$Eu_{X}(0)-Eu_{X^g}(0)-\sum_{j=1}^{r}m_{b_j}(Eu_X(b_j)-Eu_{X^g}(b_j))=(-1)^{d-1}m.$
\end{center}

\end{corollary}

\noindent\textbf{Proof.} Since $l$ is a generic linear form over $X$ and $\Sigma_{\mathcal{V}} g$ is one-dimensional, by Lemma \ref{generic linear form lemma}, $\Sigma_{\mathcal{V}} g\cap \{l=0\}=\{0\}.$ Notice that, since $l$ is generic, the local degree $m_{l,b_j}$ of $l|_{b_j}$ at the origin is precisely the multiplicity of the branch $b_j$ at the origin, which we will denote by $m_{b_j}.$
Since $B_{l,X}(0)=Eu_{X}(0)$ and $B_{l,X^g}(0)=Eu_{X^g}(0)$, we have the formula.\fim


Applying Corollary \ref{Brasselet number 1} to this case, we obtain the following consequence. 

\begin{corollary}\label{generalization 6.6}
If $H=l^{-1}(0),$ we have that \begin{center}
$B_{g,X\cap H}(0)=Eu_{X^g}(0)-\sum_{j=1}^{r}m_{b_j}(Eu_{X^g}(b_j)-B_{g,X\cap l^{-1}(\delta)}(b_j)).$ \end{center}
\end{corollary}

\noindent\textbf{Proof.} Applying Corollary \ref{Brasselet number 1} to $f=l,$ since $B_{l,X^g}(0)=Eu_{X^g}(0)$ we have the formula.\fim


\begin{observation}\label{corolario 6.6 obs}
If $l$ is a generic linear form over $\mathbb{C}^n$,  $l^{-1}(\delta)$ intersects $X\cap\{g=0\}$ transversely  and then
\begin{center}
$Eu_{X^g}(b_j)=Eu_{X^g\cap l^{-1}(\delta)}(b_j\cap l^{-1}(\delta))=B_{g,X\cap l^{-1}(\delta)\cap L}(b_j),$
\end{center} 
where the last equality is justify by Corollary 6.6 of \cite{DG} and $L$ is a generic hyperplane in $\mathbb{C}^n$ passing 
through $x_l\in l^{-1}(\delta)\cap b_j, j\in\{1,\ldots,r\}$ and $ l\in\{i_1,\ldots, i_{k(j)}\}.$ 

Denoting $B_{g,X\cap l^{-1}(\delta)\cap L}(b_j)$ by $B^{\prime}_{g,X\cap l^{-1}(\delta)}(b_j),$ the formula obtained in \ref{generalization 6.6} can be written as
\begin{center}
$B_{g,X\cap H}(0)=Eu_{X^g}(0)-\sum_{j=1}^{r}m_{b_j}(B^{\prime}_{g,X\cap l^{-1}(\delta)}(b_j)-B_{g,X\cap l^{-1}(\delta)}(b_j)).$
\end{center} 
\end{observation}

This result allows us to compare the Brasselet number $B_{g,X\cap H}(0)$ and the Euler obstruction $Eu_{X^g}(0)$ in terms of the dimension of the analytic complex space $(X,0).$

Let $I_0(X^f,\Gamma^0_{f|_X})$ be the intersection multiplicity of $X^f$ and $\Gamma^0_{f|_X},$ where $\Gamma^0_{f|_X}$ is the general relative polar curve of $f$ (see \cite{LT}). By Corollary \ref{corollary 5.2 of DG}, if $d=\dim(X),$ then\begin{eqnarray}\label{5.2 formula}
B_{f,X}(0)-B_{f,X\cap H}(0)=(-1)^{d-1}I_0(X^f,\Gamma^0_{f|_X}).
\end{eqnarray}

\begin{corollary}\label{inequalities Euler obstruction and Brasselet number}
If $H=l^{-1}(0)$ we have:
\begin{enumerate}
\item {If $d$ is even, $B_{g,X\cap H}(0)\geq Eu_{X^g}(0);$}

\item {If $d$ is odd, $B_{g,X\cap H}(0)\leq Eu_{X^g}(0).$}
\end{enumerate}
\end{corollary}
\noindent\textbf{Proof.} For $0<|\delta|\ll\epsilon\ll1$ and a generic linear form $l$ defined over $X,$ we apply  Formula (\ref{5.2 formula}) to the space $X\cap l^{-1}(\delta)$, which dimension is $d-1.$ We fix a point $x_{l}\in l^{-1}(\delta)\cap b_j,$ for each $j\in\{1,\ldots,r\}, l\in\{i_1,\ldots, i_{k(j)}\}$ and calculate the difference of Brasselet numbers around the singular point $x_{l}:$ 
 \begin{eqnarray*}
B_{g,X\cap l^{-1}(\delta)}(x_{l})-B^{\prime}_{g,X\cap l^{-1}(\delta)}(x_{l})=(-1)^{d-2}I_{x_{l}}((X\cap l^{-1}(\delta))^g,\Gamma^{x_{l}}_{g|_{X\cap l^{-1}(\delta)}}).
\end{eqnarray*}
\fim
\vspace{0.1 cm}

If $g:\mathbb{C}^n\rightarrow\mathbb{C}$ has isolated singularity at the origin, Lê and Teissier proved, in \cite{LT}, that, for $0< |\alpha|\ll\epsilon\ll 1,$ $Eu_{X^g}(0)=\chi(g^{-1}(\alpha)\cap H\cap B_{\epsilon}),$ where $X=\mathbb{C}^n$, $H$ is a generic hyperplane and $\alpha$ is a regular value of $g.$ The next result is a generalization of this result to our setting. 

Let $l$ be a generic linear form over $\mathbb{C}^n,$ $\{\mathbb{C}^n\setminus\{0\},\{0\}\}$ a Whitney stratification of $\mathbb{C}^n$ and $\{\mathbb{C}^n\setminus\{l=0\},\{l=0\},\{0\}\}$ be the good stratification of $\mathbb{C}^n$ induced by $l.$
Consider a point $x_{l}\in \{l=\delta\}\cap b_j,$ for each $j\in\{1,\ldots,r\}, l\in\{i_1,\ldots,i_{k(j)}\}$ and let $D_{x_l}$ be the closed ball with center at $x_{l}$ and radius $r_{l}, 0<|\alpha|\ll|\delta|\ll r_{l}\ll\epsilon\ll1,$ sufficiently small such that the balls $D_{x_l}$ are pairwise disjoint and the union of balls $D_j=D_{x_{i_1}}\cup\ldots\cup D_{x_{i_{k(j)}}}$ is contained in $B_{\epsilon}.$ 

\begin{corollary}
Let $H=l^{-1}(0)$ be a generic hyperplane through the origin. For $x_{l}\in \{l=\delta\}\cap b_j,$  $j\in\{1,\ldots,r\}, l\in\{i_1,\ldots,i_{k(j)}\},$ chosen as before,
\begin{center}
$Eu_{\{g=0\}}(0)=\chi(g^{-1}(\alpha)\cap H\cap B_{\epsilon})+\sum_{j=1}^{r}(-1)^{n-1}m_{b_j}(\mu(g|_{l^{-1}(\delta)},x_{l})+\mu^{\prime}(g|_{l^{-1}(\delta)},x_{l})).$
\end{center}

\end{corollary}
\noindent\textbf{Proof.} Applying Remark \ref{corolario 6.6 obs}, we obtain 

\begin{center}
$B_{g,H}(0)=Eu_{\{g=0\}}(0)-\sum_{j=1}^{r}m_{b_j}(B_{g,l^{-1}(\delta)\cap L}(x_{l})-B_{g,l^{-1}(\delta)}(x_{l})),$
\end{center} 
where $L$ is a generic hyperplane in $\mathbb{C}^n$ passing through $x_{l}.$ 
Using the definition of the Brasselet number, we have 
\begin{eqnarray*}
B_{g,H}(0)&=&\chi(H\cap g^{-1}(\alpha)\cap B_{\epsilon})\\
B_{g,l^{-1}(\delta)}(x_{l})&=&\chi(g^{-1}(\alpha)\cap l^{-1}(\delta)\cap D_{{x_l}})=1+(-1)^{n-2}\mu(g|_{l^{-1}(\delta)},x_{l})\\
B_{g,l^{-1}(t_0)\cap L}(x_{l})&=& \chi(g^{-1}(\alpha)\cap l^{-1}(\delta)\cap L\cap D_{{x_l}})=1+(-1)^{n-3}\mu(g|_{l^{-1}(\delta)\cap L},x_{l})\\&=&1+(-1)^{n-3}\mu^{\prime}(g|_{l^{-1}(\delta)},x_{l}).
\end{eqnarray*} \fim

\begin{observation}
Using the previous corollary, we can compare the difference\linebreak $Eu_{X^g}(0)-\chi(g^{-1}(\alpha)\cap H\cap B_{\epsilon})$ using the dimension of $\mathbb{C}^n.$

\begin{enumerate}
\item {If $n$ is even, $Eu_{X^g}(0)\leq\chi(g^{-1}(\alpha)\cap H\cap B_{\epsilon}).$}

\item {If $n$ is odd, $Eu_{X^g}(0)\geq\chi(g^{-1}(\alpha)\cap H\cap B_{\epsilon}).$}
\end{enumerate}
\end{observation}

In Corollary \ref{6.7 of DG}, the authors showed that if $g$ has an isolated singularity, $l$ is a generic linear form and $V_q$ is the top stratum of the stratification of $X,$ \begin{center}
$\mu^g(\Gamma_{g,l}(V_q))-\mu^l(\Gamma_{g,l}(V_q))=(-1)^d Eu_{g,X}(0)=(-1)^{d-1}(B_{g,X}(0)-Eu_X(0)).$
\end{center}

In the following, we use Corollary \ref{generalization 6.5} to present a generalization of this result to a function-germ $g:X\rightarrow\mathbb{C}$ with one-dimensional critical locus.

\begin{corollary}
We have \begin{center}
$\mu^g(\Gamma_{g,l}(V_q))-\mu^l(\Gamma_{g,l}(V_{q}))=(-1)^{d-1}\left(B_{g,X}(0)-Eu_X(0)+\sum_{j=1}^{r}m_{b_j}(Eu_X(b_j)-B_{g,X\cap\{l=\delta\}}(b_j))\right).$
\end{center}
\end{corollary}

\noindent\textbf{Proof.} Applying Corollary \ref{generalization 6.5} to $f=l$, we obtain \begin{equation}
B_{g,X}(0)-B_{l,X}(0)=(-1)^{d-1}(n_{reg}-m_{reg})-\sum_{j=1}^{r}m_{b_j}(Eu_X(b_j)-B_{g,X\cap\{l=\delta\}}(b_j)).\label{F2}
\end{equation}

By Proposition 1.14 of \cite{Ms1}, $l$ is decent with respect to $\mathcal{V}$ relative to $g$ and $g$ is decent with respect to $\mathcal{V}$ relative to $l.$ Then, we can replace $n_{reg}$ with $\mu^g(\Gamma_{g,l}(V_q))$ and $m_{reg}$ with $\mu^l(\Gamma_{g,l}(V_{q})).$ Hence we have the formula, since $B_{l,X}(0)=Eu_X(0).$ \fim

\vspace{0,5cm}

Using the last two results, we obtain another way to calculate the Brasselet number $B_{g,X}(0).$ 

\begin{proposition}\label{Generalization of 6.6 2}
For $0<|\delta|\ll\epsilon\ll 1,$\begin{center}
$B_{g,X}(0)=(-1)^{d-1}n_{reg}+Eu_{X^g}(0)-\sum_{j=1}^{r}m_{b_j}(Eu_{X^g}(b_j)-B_{g,X\cap\{l=\delta\}}(b_j)),$
\end{center}
\noindent where $n_{reg}$ is the number of stratified Morse critical points of the Morsification of \linebreak$l:X\cap g^{-1}(\delta)\cap B_{\epsilon}\rightarrow\mathbb{C}$ appearing on $X_{reg}\cap g^{-1}(\delta)\cap \{l\neq 0\}\cap B_{\epsilon}.$
\end{proposition}

\noindent\textbf{Proof.} Applying Corollary \ref{generalization 6.5} to the case where $f$ is the generic linear form $l$ and $B_{l,X}(0)=Eu_X(0),$ \begin{equation}\label{aux1}
B_{g,X}(0)=(-1)^{d-1}(n_{reg}-m_{reg})+Eu_X(0)-\sum_{j=1}^{r}m_{b_j}(Eu_X(b_j)-B_{g,X\cap\{l=\delta\}}(b_j)).
\end{equation} 

But, by Corollary \ref{Difference Euler obstructions}, \begin{equation}\label{aux2}
Eu_{X}(0)-\sum_{j=1}^{r}m_{b_j}Eu_X(b_j)=(-1)^{d-1}m_{reg}+Eu_{X^g}(0)-\sum_{j=1}^{r}m_{b_j}Eu_{X^g}(b_j).
\end{equation}

So, using equations (\ref{aux1}) and (\ref{aux2}), we obtain the formula.\fim

\begin{corollary}\label{Generalization of 6.6 3}
For $0<|\delta|\ll\epsilon\ll 1,$\begin{center}
$B_{g,X}(0)=(-1)^{d-1}n_{reg}+Eu_{X^g}(0)-\sum_{j=1}^{r}m_{b_j}(B^{\prime}_{g,X\cap l^{-1}(\delta)}(b_j)-B_{g,X\cap l^{-1}(\delta)}(b_j)),$
\end{center}
\noindent where $n_{reg}$ is the number of stratified Morse critical points of the Morsification of \linebreak$l:X\cap g^{-1}(\delta)\cap B_{\epsilon}\rightarrow\mathbb{C}$ appearing on $X_{reg}\cap g^{-1}(\delta)\cap \{l\neq 0\}\cap B_{\epsilon}.$  
\end{corollary}

\noindent\textbf{Proof.} We have the formula since, by Remark \ref{corolario 6.6 obs}, $Eu_{X^g}(b_j)=B^{\prime}_{g,X\cap l^{-1}(\delta)
}(b_j).$\fim

\begin{observation}
By Theorem \ref{4.4 DG}, $(-1)^{d-1}n_{reg}=B_{g,X}(0)-B_{g,X\cap H}(0),$ where  $H=l^{-1}(0).$ Using this equality in Corollary \ref{Generalization of 6.6 3}, we obtain \begin{center}
$B_{g,X\cap H}(0)=Eu_{X^g}(0)-\sum_{j=1}^{r}m_{b_j}(B^{\prime}_{g,X\cap l^{-1}(\delta)}(b_j)-B_{g,X\cap l^{-1}(\delta)}(b_j)),$
\end{center} 
\noindent which is the formula obtained in Remark \ref{corolario 6.6 obs}.
\end{observation}



\vspace{2cm}
(Hellen Monção de Carvalho Santana) Universidade de São Paulo, Instituto de Ciências Matemáticas e de
Computação - USP, Avenida Trabalhador São-Carlense, 400 - Centro, São Carlos, Brazil.
\it{E-mail address}: hellenmcarvalho@hotmail.com
 \end{document}